# NUMBER FIELDS WITH DISCRIMINANT $\pm 2^a 3^b$: EXAMPLES FROM THREE POINT COVERS

DAVID P. ROBERTS

Let $S$ be a finite set of primes and let $G$ be a finite group. Let $NF(S, G)$ be the set of Galois number fields $K \subset \mathbf{C}$ with $\mathrm{Gal}(K/\mathbf{Q}) \cong G$ and discriminant divisible only by primes in $S$. The sets $NF(S, G)$ are finite, by a classical theorem of Hermite.

We are interested in the following inverse Galois problem. *For given $S$ and $G$, find defining polynomials $f(x) \in \mathbf{Q}[x]$ for as many fields in $NF(S, G)$ as possible.* This problem is most interesting when $|NF(S, G)|$ can be expected to be small, so that one can reasonably aim for complete lists. Thus we are interested in "collecting number fields", the most valuable specimens being those which are least ramified, for a given Galois group $G$.

In this paper we consider only cases with $S = \{2, 3\}$. If $G$ has a faithful permutation representation of degree $\leq 7$, we have completely identified $NF(S, G)$ previously with Jones [JR1], [Jon]. Here we supplement these complete lists with incomplete lists of fields with larger groups $G$, involving the simple groups $PSL_2(9)$, $SL_2(8)$, $SL_3(3)$, $A_8$, $PSp_4(3)$, $A_9$, $Sp_6(2)$, and $A_{32}$.

We construct all our fields by the standard technique of specializing three point covers. However our requirement that the final field discriminant be of the form $\pm 2^a 3^b$ is extremely restrictive. In the main, it restricts consideration to three point covers with bad reduction at 2 and 3 only. To keep ramification within $\{2, 3\}$, the specialization point $\tau$ then has to be chosen judiciously as well. Throughout we used *Mathematica*, supplemented by the *discf*, *factorpadic*, and *polred* functions of *Pari*.

Of the three point covers which we use, for ten it was hard to find a defining polynomial. The polynomials $f_{9a}$, $f_{9b}$, $f_{13a}$, and $f_{27a}$ were found by Elkies, Matzat, Malle, and Häfner, respectively, [Elk1], [Mat], [Mal1], [Häf]. The polynomials $f_{13b}$, $f_{13c}$, $f_{27b}$, $f_{27c}$, $f_{27d}$, and $f_{28}$ are new here. We do not go into much detail about how we found the new polynomials, our methods being standard. Rather, we emphasize various aspects of the specialization step.

Section 1 provides an overview of the fields constructed in this paper. Section 2 and 3 discuss three point covers and especially specialization. Sections 4-10 discuss the individual constructions, sorted approximately by $|G|$; each of these sections also informally discusses a more theoretical issue, as indicated by the section titles. Section 11 indicates a direction for further research we feel is important and promising. We conclude in Section 12 with two speculative conjectures on the general nature of the numbers $|NF(S, G)|$.

Throughout we work as much as possible in a purely number-theoretic setting. However the reader should realize that a systematic study of the sets $NF(S, G)$ would fundamentally involve both the theory of motives and the theory of automorphic representations:

*Motives.* Let $X$ be a variety over $\mathbf{Q}$ with bad reduction within a set $S_0$. Let $\ell$ be a prime number and put $\mathbf{F}_\ell = \mathbf{Z}/\ell$. Then $\mathrm{Gal}(\overline{\mathbf{Q}}/\mathbf{Q})$ acts on the set $H^w(X(\mathbf{C}), \mathbf{F}_\ell)$ with kernel say $U$. Then any Galois subfield $K \subset \mathbf{C}$ fixed by $U$ is ramified within $S_0 \cup \{\ell\}$. In fact, most of the new number fields here arise in this way as mod 2 and/or mod 3 Galois representations associated to pure motives defined by Katz [Kat]; we explain this connection in [Rob2].

*Automorphic Representations.* When $G$ is solvable, the sets $NF(S, G)$ are in principle accessible by abelian class field theory. Extending Serre's $GL_2$ conjecture [Ser2], Gross has a conjectural Langlands-type non-abelian class field theory for quite general groups, [Gro1], [Gro2]. It predicts in particular that for a given prime $p$, $NF(\{p\}, G)$ is non-empty for suitable non-solvable $G$. For example, in [Gro2] it is predicted that $NF(\{5\}, G_2(5))$ is non-empty.

The motivic approach to $NF(S, G)$ yields defining equations only in very special cases, such as the classical case of taking $\ell$-torsion points on elliptic curves. The automorphic approach to $NF(S, G)$ never yields defining equations. It is for these reasons that we emphasize explicit defining equations here.

## Contents

## 1. Summary of fields constructed

To keep our investigation of manageable size, we focus attention on groups $G$ of the form $H.A$, where $H$ is non-abelian simple and the natural map $A \to \mathrm{Out}(H)$ is injective. In passing, however, we construct interesting fields for $G$ not in this class. For example, in Sections 4-7, some interesting solvable fields appear. Also, in Section 5 some $G$ involving two copies of $SL_2(8)$ appear, and many other such $H^n.A$ could be constructed just by specializing at non-rational points $\tau$. Also, the cover given by $f_{27d}(t, x)$ in Section 7 has natural lifts to covers with Galois group $GSp_4(3) = 2.PSp_4(3).2$.

Table 1.1 lists some groups $G$ of the form $H.A$. Groups $G$ are placed in the same block iff they have the same $H$. We give a characteristic 0 description, a characteristic 2 description, and/or a characteristic 3 description. The Atlas [Atlas] sometimes gives even more descriptions, e.g. $SL_3(2).2 \cong PGL_2(7)$. The column $N$ gives the degree of a minimal faithful permutation representation.

An entry under $\#_G$ of the form $\bullet x$ summarizes results from [JR1] and/or [Jon]. In this case, $|NF(\{2,3\},G)| = x$. The complete search for $GL_3(2)/A_7/S_7$ fields took three months, but an analogous complete search for $GL_3(2).2/A_8/S_8$ fields would take somewhat more than a million years [JR2]. This is why we are shifting attention to non-exhaustive but still systematic ways of constructing fields in a given $NF(S,G)$.

Table 1.1. Lower bounds for $|NF(\{2,3\},G)|$.

| $\lvert G\rvert$ | 0 | 2 | 3 | $N$ | $\#_G$ | -6 | -3 | -2 | -1 | 2 | 3 | 6 |
|---|---|---|---|---|---|---|---|---|---|---|---|---|
| $60 = 2^2 3 5$ | $A_5$ | $SL_2(4)$ | | 5 | $\bullet 0$ | | | | | | | |
| $120 = 2^3 3 5$ | $S_5$ | $\Sigma L_2(4)$ | | 5 | $\bullet 5$ | 1 | 0 | 0 | 0 | 2 | 1 | 1 |
| $168 = 2^3 3 7$ | | $SL_3(2)$ | | 7 | $\bullet 0$ | | | | | | | |
| $336 = 2^4 3 7$ | | $SL_3(2).2$ | | 8 | | | | | | | | |
| $360 = 2^3 3^2 5$ | $A_6$ | $Sp_4(2)'$ | $PSL_2(9)$ | 6 | $\bullet 4$ | | | | | | | |
| $720 = 2^4 3^2 5$ | $S_6$ | $Sp_4(2)$ | $PGO_4^-(3)$ | 6 | $\bullet 27$ | 1 | 1 | 0 | 2 | 13 | 4 | 6 |
| $720 = 2^4 3^2 5$ | | | $PGL_2(9)$ | 10 | | | | | | | | |
| $720 = 2^4 3^2 5$ | $M_{10}$ | | | 10 | 10 | 0 | 0 | 0 | 0 | 10 | 0 | 0 |
| $1,440 = 2^5 3^2 5$ | | | $P\Gamma L_2(9)$ | 10 | 60 | | | | | | | |
| $504 = 2^6 3^2 7$ | | $SL_2(8)$ | | 9 | 3 | | | | | | | |
| $1,512 = 2^6 3^3 7$ | | $\Sigma L_2(8)$ | | 9 | 63 | | | | | | | |
| $2,520 = 2^3 3^2 5 7$ | $A_7$ | | | 7 | $\bullet 0$ | | | | | | | |
| $5,040 = 2^4 3^2 5 7$ | $S_7$ | | | 7 | $\bullet 10$ | 4 | 1 | 1 | 4 | 0 | 0 | 0 |
| $5,616 = 2^4 3^3 13$ | | | $SL_3(3)$ | 13 | 6 | | | | | | | |
| $11,232 = 2^5 3^3 13$ | | | $SL_3(3).2$ | 26 | 85 | 6 | 6 | 64 | 6 | 1 | 0 | 2 |
| $6,048 = 2^5 3^3 7$ | | $G_2(2)'$ | $SU_3(3)$ | 28 | | | | | | | | |
| $12,096 = 2^6 3^3 7$ | | $G_2(2)$ | $SU_3(3).2$ | 28 | | | | | | | | |
| $20,160 = 2^6 3^2 5 7$ | $A_8$ | $GL_4(2)$ | | 8 | | | | | | | | |
| $40,320 = 2^7 3^2 5 7$ | $S_8$ | $SO_6^+(2)$ | | 8 | 2 | 2 | 0 | 0 | 0 | 0 | 0 | 0 |
| $25,920 = 2^6 3^4 5$ | $W(E_6)'$ | $SU_4(2)$ | $PSp_4(3)$ | 27 | 21 | | | | | | | |
| $51,840 = 2^7 3^4 5$ | $W(E_6)$ | $SO_6^-(2)$ | $SO_5(3)$ | 27 | 124 | 18 | 54 | 30 | 13 | 3 | 0 | 6 |
| $181,440 = 2^6 3^4 5 7$ | $A_9$ | | | 9 | 10 | | | | | | | |
| $362,880 = 2^7 3^4 5 7$ | $S_9$ | | | 9 | 25 | 0 | 2 | 7 | 1 | 5 | 6 | 4 |
| $1,451,520 = 2^9 3^4 5 7$ | $W(E_7)'$ | $Sp_6(2)$ | | 28 | 34 | | | | | | | |
| $2.63 \times 10^{35} \approx 32!$ | $S_{32}$ | | | 32 | 1 | 1 | 0 | 0 | 0 | 0 | 0 | 0 |

The main focus of the present paper is the following theorem.

**Theorem 1.1.** *For $G = M_{10}$, $P\Gamma L_2(9)$, $SL_2(8)$, $\Sigma L_2(8)$, $S_8$, $PSp_4(3)$, $SO_5(3)$, $A_9$, $S_9$, $Sp_6(2)$, and $S_{32}$, the polynomials presented in Section 4-10 give $\#_G$ elements of $NF(\{2,3\},G)$, with $\#_G$ as on Table 1.1.*

In general, suppose given separable polynomials $f_i \in \mathbf{Q}[x]$, which one expects have distinct splitting fields $K_i$ in some $NF(S,G)$. After a variable change, one can assume all the $f_i$ are monic polynomials in $\mathbf{Z}[x]$. To prove that the $K_i$ are in $NF(S,G)$ and distinct one has to do three things:

*1. Verify that the ramification set $S_i$ of $K_i$ is really in $S$.* Let $L_i = \mathbf{Q}[x]/f_i(x)$. Let $D_i$ be the polynomial discriminant $f_i$ and let $d_i$ be the field discriminant of $L_i$. One has to show that the $p \notin S$ dividing $D_i$ do not divide $d_i$. This may require substantial computation in general, but in our case we do it uniformly without computation by the standard fact 2.4.

*2. Verify that the Galois group* $G_i := \mathrm{Gal}(K_i/\mathbf{Q})$ *is indeed isomorphic to* $G$. One can very easily compute lower bounds by means of Frobenius elements; we briefly indicate how this goes in Section 4, and then give no more details. On a heuristic level, one can expect that the lower bound obtained after substantial computation is always exact. On a rigorous level, to get upper bounds may require very substantial computation, like those we carry out for related purposes in Sections 7 and 9. In our case, we get the needed upper bounds uniformly without computation from the standard fact 2.2.

*3. Verify that the* $K_i$ *are distinct.* Often the group $G$ has up to isomorphism exactly one faithful permutation representation of a given degree $N$, and all given defining polynomials $f_i$ have degree $N$. In this case, for each pair $i \neq i'$ one needs to find a prime $p$, not dividing $D_i D_{i'}$, for which the factorization partitions of $f_i$ and $f_{i'}$ over $\mathbf{F}_p$ are distinct. In general the situation may be very slightly more complicated, for example by sextic twinning in Section 4 and projective twinning in Section 6. But even here, one needs to just find $p$ such that the factorization partitions come from elements of different orders in the symmetric group $S_N$. We have done this, but do not present any details here.

In the rest of the paper we complete the statement of Theorem 1.1 by presenting the defining equations. The proof of Theorem 1.1 is completed simultaneously, since, as we have just explained, the theorem is essentially self-proving from its full statement. Of course, it would be pointless to put on blinders, restricting ourselves only to establishing those facts literally contributing to Theorem 1.1. Rather, we present a fuller picture of both our particular examples and the general technique of constructing fields in a given $NF(S, G)$ by means of specializing three point covers.

For $K \in NF(\{2,3\}, H.A)$ let $K^s$ be the maximal solvable subfield, so that $\mathrm{Gal}(K^s/\mathbf{Q}) \cong A$. When $|A| = 2$, $K^s = \mathbf{Q}(\sqrt{d})$ for some $d \in \{-6, -3, -2, -1, 2, 3, 6\}$. The last block of columns sorts the fields considered here according to this "discriminant class" $d$. One place where the nature of $A$ plays an important role is in comparing number fields with automorphic forms. For example non-totally-real $PGL_2(9)$ and $SL_2(8)$ fields should come from classical modular forms, according to Serre's conjecture, while $P\Gamma L_2(9)$ and $\Sigma L_2(8)$ fields certainly do not. On the other hand, one would expect that these $P\Gamma L_2(9)$ and $\Sigma L_2(8)$ fields come from Hilbert modular forms. Somewhat similarly, Gross has explained to me that $PGL_n(\ell)$ and $PGL_n(\ell).2$ fields with $d > 0$ should be relatively rare, while $PGL_n(\ell).2$ fields with $d < 0$ should be more common; here $n \geq 3$, like in our case $PGL_3(3) = SL_3(3)$.

On Table 1.1, there are also some blanks under $\#_G$. In these cases we know of no defining equations for corresponding fields; these rows are included as specific suggestions for further research. From looking at mod $\ell$ Galois representations of various curves, we know that there are more fields for $PGL_2(9)$, $SL_2(8)$, $SU_3(3).2$ and $SO_5(3)$ than listed here. From modest computer searches we also know that there are more fields for $S_8$ than listed here.

In [JR1] we analyzed the local behavior at 2 and 3 of low degree fields in complete detail. In particular, all $5 + 4 + 27$ fields on the top lines of Table 1.1 are wildly ramified at both 2 and 3, mostly very wildly ramified, as explained in Section 3.3 there. Similarly, all 10 $S_7$ fields are wildly ramified at both 2 and 3. In contrast, some of the degree $\leq 13$ fields here are wildly ramified at only one of these two primes; see 4.1, 5.1, 5.2. We have not investigated ramification at 2 and 3 in the degree 26-28 fields.

## 2. Three point covers and specialization

Let $F$ be a field and work with smooth projective curves over the fixed projective line $\mathbf{P}^1_F = \mathrm{Spec}(F[t]) \cup \{\infty\}$. A three point cover over $F$ is a finite separable cover $X \to \mathbf{P}^1_F$, ramified only above the three points 0, 1, and $\infty$. A defining polynomial for $X$ is a polynomial $f(t,x)$ with $F(X) = F(t)[x]/f(t,x)$. An important issue extensively addressed in the literature is how one explicitly constructs defining equations; throughout this paper we take the defining equations simply as given, and the bulk of this section is devoted to setting up notation. Standard references for three point covers include the books [Mat], [SD], [Sch], and [Völ].

We need first some notation with respect to the universal base curve. Let $T = \mathrm{Spec}(\mathbf{Z}[t, 1/t, 1/(t-1)])$ be the thrice-punctured projective line, so that for any field $F$, $T(F) = F - \{0,1\}$. The topological space $T(\mathbf{C})$ plays a central role. We work with the fundamental group

$$\pi_1 := \pi_1(T(\mathbf{C}), 1/2) = \langle \gamma_0, \gamma_1, \gamma_\infty | \gamma_0\gamma_1\gamma_\infty = 1\rangle.$$

Here $\gamma_0$ comes from the counterclockwise circle about 0 with radius 1/2; $\gamma_1$ comes from the counterclockwise circle about 1 with radius 1/2. The group $\pi_1$ is free on the two generators $\gamma_0$ and $\gamma_1$; the third element $\gamma_\infty := (\gamma_0\gamma_1)^{-1}$ is introduced so as to treat the three cusps in the same way.

To keep control of groups even when specializing, it is convenient to remove branch-cuts as follows.

$$T(\mathbf{C})^{\mathrm{cut}} = \mathbf{C} - \{1 + t(1 - i\pi)\}_{t \in \mathbf{R}_{\geq 0}} - \{t(-1 - i\pi)\}_{t \in \mathbf{R}_{\geq 0}}.$$

The cut plane $T(\mathbf{C})^{\mathrm{cut}}$ is simply connected. We have chosen the branch-cuts so that all of $\overline{\mathbf{Q}} - \{0,1\}$ is in $T(\mathbf{C})^{\mathrm{cut}}$.

A simple partition of $T(\mathbf{Q}_p)$ into subregions plays a prominent role in the analysis of ramification. For $\tau \in T(\mathbf{Q}_p)$ write $\tau = -A/C$ with $A, C \in \mathbf{Z}_p$ not both divisible by $p$. Define $B \in \mathbf{Z}_p$ by $A + B + C = 0$. Define

$$\begin{aligned} \mathrm{ord}_0(\tau) &:= \mathrm{ord}_p(A) \\ \mathrm{ord}_1(\tau) &:= \mathrm{ord}_p(B) \\ \mathrm{ord}_\infty(\tau) &:= \mathrm{ord}_p(C). \end{aligned}$$

The partition is

(2.1)
$$T(\mathbf{Q}_p) = T(\mathbf{Q}_p)^{\mathrm{gen}} \coprod \left(\coprod_{i=1}^{\infty} T(\mathbf{Q}_p)^{0,i}\right) \coprod \left(\coprod_{i=1}^{\infty} T(\mathbf{Q}_p)^{1,i}\right) \coprod \left(\coprod_{i=1}^{\infty} T(\mathbf{Q}_p)^{\infty,i}\right).$$

Here $\mathrm{ord}_c$ takes the value $i$ exactly on $T(\mathbf{Q}_p)^{c,i}$; so for $c$ fixed and $i$ increasing the $T(\mathbf{Q}_p)^{c,i}$ are smaller and smaller annuli, all centered at the cusp $c$. The generic piece, empty for $p = 2$, is where all three $\mathrm{ord}_c$ take the value 0.

For the rest of this section a degree $N$ three point cover $X \to \mathbf{P}^1_F$ is fixed. A particular defining polynomial $f(t,x)$ is fixed too. Objects $\Delta$, $a$, $b$, $c(t)$, and $\Sigma$ below depend on $f$. Otherwise the objects $M$, $G$, ... depend only on $X$. We will use this notation systematically in the sequel, with indices when discussing particular examples. For example, once we have denoted a particular cover $X_{27d}$, then $M_{27d}$ is automatically its monodromy group.

The polynomial discriminant of $f(t,x)$, with respect to the variable $x$, factors canonically as

$$D(t) = \Delta t^a (t-1)^b c(t)^2.$$

with $c(t) \in F[t]$ monic and prime to $t(t-1)$. Here $\Delta \in F^\times$ and $a$ and $b$ are integers. Let $\Sigma$ be the subvariety of $\mathbf{P}^1_F$ corresponding to the roots of $c(t)$. Always there is a different defining polynomial $f^*(t,x)$ with $c^*(t)$ relatively prime to $c(t)$. So the factor $c(t)^2$ plays essentially no role in our situation; actually $c(t)$ is identically 1 in most of our examples.

For $\tau$ in $F - \{0,1\}$ let $X_\tau = \mathrm{Spec}(L_\tau)$ be the corresponding fiber. So $L_\tau = F(X_\tau)$ is a finite separable algebra over $F$. For $\tau \notin \Sigma(F)$, $L_\tau = F[x]/f(\tau,x)$; so in this case $L_\tau$ is a field iff $f(\tau,x)$ is irreducible.

Suppose henceforth that $F \subseteq \mathbf{C}$. Let $X(\mathbf{C})^{\mathrm{cut}}$ be the inverse image of $T(\mathbf{C})^{\mathrm{cut}}$ in $X(\mathbf{C})$. Let $C$ be set of components of $X(\mathbf{C})^{\mathrm{cut}}$, i.e. $C = \pi_0(X(\mathbf{C})^{\mathrm{cut}})$. The group $\pi_1$ acts naturally on $C$. The image of $\pi_1$ in the symmetric group $\mathrm{Sym}(C)$ is called the monodromy group $M$. For $c \in \{0,1,\infty\}$ the image of $\gamma_c$ in $M$ is denoted $m_c$; its order is denoted $e_c$.

Let $\overline{\mathbf{Q}}$ be the algebraic closure of $\mathbf{Q}$ in $\mathbf{C}$. Suppose henceforth that $F \subseteq \overline{\mathbf{Q}}$; this is the essential case anyway. For $\tau \in F - \{0,1\}$, the branch cuts let us canonically identify the fiber $X_\tau(\mathbf{C})$ with the component set $C$; we use this identification without comment in the sequel. But also $X_\tau(\mathbf{C})$ is identified with the set of homomorphisms $\mathrm{Hom}(L_\tau, \overline{\mathbf{Q}})$; for $\tau \notin \Sigma(F)$, $X_\tau(\mathbf{C})$ is identified with the set of roots $X_\tau(\mathbf{C}) \subset \overline{\mathbf{Q}}$ of $f(\tau,x)$. Either way, one sees a natural action of $\mathrm{Gal}(\overline{\mathbf{Q}}/F)$ on $X_\tau(\mathbf{C}) = X_\tau(\overline{\mathbf{Q}})$. The image of $\mathrm{Gal}(\overline{F}/F)$ in the symmetric group $\mathrm{Sym}(C)$ is the Galois group $G_\tau$.

Let $M^*$ be the normalizer of $M$ in the symmetric group $\mathrm{Sym}(C)$. A basic fact is that

(2.2) each $G_\tau$ is contained in $M^*$.

Another basic fact is that all the biggest $G_\tau$ coincide, and we call this common group the Galois group $G$.

Suppose now that $F \subset \overline{\mathbf{Q}}$ is a number field of degree $d$. Then one has $d$ conjugate covers $X^\delta \to F^\delta$. Restricting scalars, one gets a single cover of degree $dN$ over the field $\mathbf{Q}$. On the level of defining polynomials, this operation of restricting scalars is the passage from $f(t,x) \in F(t,x)$ to $\prod f^\delta(t,x) \in \mathbf{Q}[t,x]$. Since our object is to construct Galois extensions of $\mathbf{Q}$, this is a natural thing for us to do; in the rest of this section, we take $F = \mathbf{Q}$.

Let $K_\tau$ be the splitting field in $\mathbf{C}$ of $L_\tau$. So if $\tau \notin \Sigma(\mathbf{Q})$, $K_\tau$ is the field generated by the roots of $f(\tau,x)$ in $\mathbf{C}$. The fields $K_\tau$ are the fields which we emphasize in the introduction and Sections 1,12; we have set things up so that $G_\tau = \mathrm{Gal}(K_\tau/\mathbf{Q})$. On the other hand, for most of this paper we focus on the degree $N$ algebras $L_\tau$ which arise naturally and are computationally more accessible. Thus, for example, we sometimes give the algebra discriminant of $L_\tau$, as computed by *discf*. Determining the discriminant of $K_\tau$ would typically require more calls to *discf*, not possible if degrees are too large, plus theoretical work, based on the conductor-discriminant formula.

Let $S$ be the set of primes at which the cover $X \to \mathbf{P}^1_{\mathbf{Q}}$ has bad reduction. Always $S$ contains all primes dividing $e_0 e_1 e_\infty$. Always $S$ is contained in the primes dividing numerator or denominator of $\Delta$. By adjusting the defining equation, one can always make $S$ exactly the primes dividing numerator or denominator of $\Delta$.

A basic fact is that all the fields $K_\tau$ are tamely ramified outside $S$. To be more precise, recall that for a Galois number field $K \subset \mathbf{C}$, tamely ramified at $p$, one has a conjugacy class $r_p \in \mathrm{Gal}(K/\mathbf{Q})$, measuring the ramification at $p$. In our case, for $p \notin S$, one has an explicit formula for the ramification class $r_{\tau,p}$ as a conjugacy class in $G \supseteq G_\tau$:

$$r_{\tau,p} = \begin{cases} 1 & \text{if } \tau \in T(\mathbf{Q}_p)^{\mathrm{gen}} \\ m_c^i & \text{if } \tau \in T(\mathbf{Q}_p)^{c,i}. \end{cases} \tag{2.3}$$

Thus $L_\tau$ is unramified at $p$ iff $\tau \in T(\mathbf{Q}_p)^{\mathrm{gen}}$ or $\tau \in T(\mathbf{Q}_p)^{c,i}$ with $m_c^i = 1$. Also

$$\text{if } \tau \in T(\mathbf{Q}_p)^{c,i} \text{ then } \mathrm{ord}_p(\mathrm{disc}(L_\tau)) = N - |C/m_c^i|, \tag{2.4}$$

the last term being the number of orbits of $m_c^i$ on the component set $C$.

## 3. Specialization lists for $S = \{2, 3\}$

Let $p$, $q$, and $r$ be positive integers or $\infty$. Let $S$ be a finite set of primes. Define $T_{p,q,r}(\mathbf{Z}^S)$ to be the set of those rational numbers $\tau \in \mathbf{Q} - \{0, 1\}$ such that

$$\tau \in T(\mathbf{Q}_p)^{\mathrm{gen}} \coprod \left( \coprod_{i=1}^{\infty} T(\mathbf{Q}_p)^{0,pi} \right) \coprod \left( \coprod_{i=1}^{\infty} T(\mathbf{Q}_p)^{1,qi} \right) \coprod \left( \coprod_{i=1}^{\infty} T(\mathbf{Q}_p)^{\infty,ri} \right)$$

for all $p$ not in $S$. Here if an index is $\infty$, then the corresponding cuspidal summand is by definition empty.

Explicitly, $\tau$ is in $T_{p,q,r}(\mathbf{Z}^S)$ iff there exist integers $(a, b, c, x, y, z)$ as follows: $a$, $b$, and $c$ are divisible only by primes in $S$; $x$, $y$, and $z$ are divisible only by primes not in $S$;

$$\begin{aligned} \tau &= -\frac{ax^p}{cz^r}\ ; \\ ax^p + by^q + cz^r &= 0. \end{aligned}$$

The numbers $(a, b, c, x, y, z)$ are then uniquely determined under the following auxiliary normalization conditions: $A = ax^p$, $B = by^q$, and $C = cz^r$ are relatively prime; $x$, $y$, and $z$ are positive; two of $A$, $B$, and $C$ are positive. In short, identifying $T_{p,q,r}(\mathbf{Z}^S)$ requires finding the solutions of a generalized Fermat equation.

The discussion at the end of Section 2 makes clear why we are interested in the sets $T_{p,q,r}(\mathbf{Z}^S)$. Namely, let $X \to \mathbf{P}^1_{\mathbf{Q}}$ be a three point cover with bad reduction set within $S$ and local monodromy orders $p = e_0$, $q = e_1$, $r = e_\infty$. Then the specialized fields $K_\tau$ have bad reduction within $S$ exactly if $\tau \in T_{p,q,r}(\mathbf{Z}^S)$. This important statement is called the Chevalley-Weil theorem for $M$-curves in [Dar].

If $\{p, q, r\} = \{2, 3, 5\}$ then the monodromy group of $M$ is the alternating group $A_5$. If otherwise $1/p + 1/q + 1/r \geq 1$ then the monodromy group $M$ is solvable. Thus we are interested here exclusively in the hyperbolic case

$$\frac{1}{p} + \frac{1}{q} + \frac{1}{r} < 1. \tag{3.1}$$

Here the sets $T_{p,q,r}(\mathbf{Z}^S)$ are known to be finite [DG], the proof appealing to Faltings' general finiteness theorem for curves.

In the application to specializing three point covers, only cases satisfying the condition

$$\text{all primes dividing } pqr \text{ are in } S \tag{3.2}$$

arise.

For $S$ fixed, put

$$T_h(\mathbf{Z}^S) = \bigcup_{(p,q,r)} T_{p,q,r}(\mathbf{Z}^S) \tag{3.3}$$

the union being over $(p, q, r)$ satisfying 3.1 and 3.2. Write $(p, q, r)|(p', q', r')$ iff $p|p'$, $q|q'$, and $r|r'$, with the convention (any positive integer)$|\infty$. Then

$$T_{p,q,r}(\mathbf{Z}^S) \supseteq T_{p',q',r'}(\mathbf{Z}^S)$$

if $(p, q, r)|(p', q', r')$. Thus 3.3 can be written as a finite union over $(p, q, r)$ minimal with respect to divisibility, and $T_h(\mathbf{Z}^S)$ is also finite.

To avoid confusion, we should mention that it is natural in other contexts to consider a larger set $T_H(\mathbf{Z}^S)$ by letting $(p, q, r)$ on the right side of 3.3 run over triples satisfying 3.1 but not necessarily 3.2. The ABC conjecture would say that each $T_H(\mathbf{Z}^S)$ is finite, but this finiteness is not known. The group $S_3 = \mathrm{Sym}(\{0, 1, \infty\})$ acts naturally on the sets $T_h(\mathbf{Z}^S)$ and $T_H(\mathbf{Z}^S)$. If $2 \in S$, one has the three-element orbit $\{-1, 1/2, 2\}$; otherwise all orbits have six elements. A case which has received a lot of attention recently is the case $S = \emptyset$. The corresponding set $T_H(\mathbf{Z})$ is known to contain ten $S_3$-orbits [Beu], and conjectured not to contain any more [DG], [Dar].

Suppose now that $S$ contains at least two primes. Then, by taking division covers associated to Katz motives, one gets covers $X$ with monodromy groups $M$ involving infinitely many finite simple groups, all with bad reduction within $S$; see [Rob2]. Let $(p, q, r)$ be the local monodromy degrees. Certainly, from hypergeometric motives only, each possible $p$ arises infinitely often. We suspect that in fact each $(p, q, r)$ satisfying 3.1 and 3.2 arises infinitely often.

The set $T_{2,3,\infty}(\mathbf{Z}^{\{2,3\}})$ has been completely identified [Cog]; it has 81 elements giving rise to 56 orbits. We have carried out a several-day computer search for elements of $T_h(\mathbf{Z}^{\{2,3\}})$, with the expectation that each $\tau$ found will be involved in the construction of infinitely many essentially distinct number fields $K_\tau$. We have found 45 more orbits. Representatives of these 101 orbits are given in the next two tables. Each representative $\tau$ is placed under $(p, q, r)$ with $\tau \in T_{p,q,r}(\mathbf{Z}^{\{2,3\}})$ and $(p, q, r)$ maximal with respect to this property.

Our computer search shows that Table 3.2 is complete with respect to solutions with $|ax^p|$, $|bx^q|$, and $|cz^r|$ all $\leq 10^9$. To get larger solutions, we permuted cusps

TABLE 3.1. The 56 $S_3$-orbits of $T_h(\mathbf{Z}^{\{2,3\}})$ needed for $T^*_{2,3,\infty}$

| $a{\cdot}x^p$ | $b{\cdot}y^q$ | $c{\cdot}1$ |
|---|---|---|
| $1$ | $1$ | $-2{\cdot}1$ |
| $2{\cdot}1$ | $1$ | $-3{\cdot}1$ |
| $3{\cdot}1$ | $1$ | $-2^2{\cdot}1$ |
| $2^3{\cdot}1$ | $1$ | $-3^2{\cdot}1$ |
| $5^3$ | $3{\cdot}1$ | $-2^7{\cdot}1$ |
| $-5^2$ | $2^4{\cdot}1$ | $3^2{\cdot}1$ |
| $-5^2$ | $2^3 3{\cdot}1$ | $1$ |
| $5^2$ | $2{\cdot}1$ | $-3^3{\cdot}1$ |
| $-7^2$ | $2^4 3{\cdot}1$ | $1$ |
| $7^2$ | $2^5{\cdot}1$ | $-3^4{\cdot}1$ |
| $2{\cdot}11^2$ | $1$ | $-3^5{\cdot}1$ |
| $-17^2$ | $2^5 3^2{\cdot}1$ | $1$ |
| $-131^2$ | $5^6$ | $2^9 3{\cdot}1$ |
| $-11^3$ | $2{\cdot}5^4$ | $3^4{\cdot}1$ |
| $2{\cdot}7^2$ | $-5^3$ | $3^3{\cdot}1$ |
| $11^2$ | $-5^3$ | $2^2{\cdot}1$ |
| $13^2$ | $-2{\cdot}5^3$ | $3^4{\cdot}1$ |
| $2^2{\cdot}5^2$ | $-7^3$ | $3^5{\cdot}1$ |
| $17^2$ | $-7^3$ | $2 \cdot 3^3{\cdot}1$ |
| $-19^2$ | $7^3$ | $2 \cdot 3^2{\cdot}1$ |
| $19^2$ | $5^3$ | $-2 \cdot 3^5{\cdot}1$ |
| $13^2$ | $7^3$ | $-2^9{\cdot}1$ |
| $3^3{\cdot}7^2$ | $-11^3$ | $2^3{\cdot}1$ |
| $37^2$ | $-2^2{\cdot}7^3$ | $3{\cdot}1$ |
| $-3{\cdot}23^2$ | $11^3$ | $2^8{\cdot}1$ |
| $35^2$ | $-13^3$ | $2^2 3^5{\cdot}1$ |
| $2^2{\cdot}23^2$ | $-13^3$ | $3^4{\cdot}1$ |
| $-47^2$ | $13^3$ | $2^2 3{\cdot}1$ |

| $a{\cdot}x^2$ | $b{\cdot}y^3$ | $c{\cdot}1$ |
|---|---|---|
| $-2^3{\cdot}17^2$ | $5^3$ | $3^7{\cdot}1$ |
| $59^2$ | $-3{\cdot}11^3$ | $2^9{\cdot}1$ |
| $61^2$ | $3{\cdot}5^3$ | $-2^{12}{\cdot}1$ |
| $-71^2$ | $17^3$ | $2^7{\cdot}1$ |
| $-2^3 3^2{\cdot}13^2$ | $23^3$ | $1$ |
| $11^2$ | $23^3$ | $-2^{12} 3{\cdot}1$ |
| $73^2$ | $23^3$ | $-2^3 3^7{\cdot}1$ |
| $143^2$ | $-3{\cdot}19^3$ | $2^7{\cdot}1$ |
| $2^3{\cdot}73^2$ | $-35^3$ | $3^5{\cdot}1$ |
| $107^2$ | $-3^2{\cdot}17^3$ | $2^{15}{\cdot}1$ |
| $-215^2$ | $19^3$ | $2 \cdot 3^9{\cdot}1$ |
| $-253^2$ | $2^9{\cdot}5^3$ | $3^2{\cdot}1$ |
| $-359^2$ | $3{\cdot}35^3$ | $2^8{\cdot}1$ |
| $545^2$ | $-2{\cdot}53^3$ | $3^6{\cdot}1$ |
| $595^2$ | $-73^3$ | $2^4 3^7{\cdot}1$ |
| $-3{\cdot}389^2$ | $2{\cdot}61^3$ | $1$ |
| $-827^2$ | $73^3$ | $2^{15} 3^2{\cdot}1$ |
| $955^2$ | $-97^3$ | $2^3 3^4{\cdot}1$ |
| $1871^2$ | $-3^2{\cdot}73^3$ | $2^9{\cdot}1$ |
| $2359^2$ | $47^3$ | $-2^5 3^{11}{\cdot}1$ |
| $2681^2$ | $-193^3$ | $2^4 3^4{\cdot}1$ |
| $-2{\cdot}2761^2$ | $239^3$ | $3^{13}{\cdot}1$ |
| $8549^2$ | $-3^5{\cdot}67^3$ | $2^3{\cdot}1$ |
| $-23053^2$ | $505^3$ | $2^{27} 3{\cdot}1$ |
| $2{\cdot}21395^2$ | $-971^3$ | $3^8{\cdot}1$ |
| $39151^2$ | $-1153^3$ | $2^5 3^5{\cdot}1$ |
| $-3{\cdot}48383^2$ | $1915^3$ | $2^{13}{\cdot}1$ |
| $-2{\cdot}184211^2$ | $4079^3$ | $3{\cdot}1$ |

and applied the following "base-change" maps iteratively.

$$\begin{aligned} f_2 : T_{m,m,n}(\mathbf{Z}^S) &\to T_{m,2,2n}(\mathbf{Z}^{S\cup\{2\}}) \\ \tau &\mapsto 4\tau(1-\tau) \\ (A,B,C) &\mapsto (4AB, (2A+C)^2, -C^2) \end{aligned}$$

$$\begin{aligned} f_3 : T_{2n,2,n}(\mathbf{Z}^S) &\to T_{3,2,2n}(\mathbf{Z}^{S\cup\{3\}}) \\ \tau &\mapsto \frac{(4\tau-1)^3}{27\,\tau} \\ (A,B,C) &\mapsto ((4A+C)^3, (8A-C)^2B, -27AC^2) \end{aligned}$$

$$\begin{aligned} f_4 : T_{3n,3,n}(\mathbf{Z}^S) &\to T_{3,2,3n}(\mathbf{Z}^{S\cup\{2\}}) \\ \tau &\mapsto \frac{(9\tau-1)^3(1-\tau)}{64\,\tau} \\ (A,B,C) &\mapsto (B(9A+C)^3, (27A^2+18AC-C^2)^2, 64AC^3) \end{aligned}$$

Table 3.2. 45 more $S_3$-orbits of $T_h(\mathbf{Z}^{\{2,3\}})$

| $a{\cdot}x^p$ | $b{\cdot}y^q$ | $c{\cdot}z^r$ |
|---|---|---|
| $2^5{\cdot}1$ | $-3{\cdot}5^3$ | $7^3$ |
| $23^2$ | $-5^4$ | $2^5 3{\cdot}1$ |
| $7^2$ | $-5^4$ | $2^6 3^2{\cdot}1$ |
| $-29^2$ | $5^4$ | $2^3 3^3{\cdot}1$ |
| $47^2$ | $-7^4$ | $2^6 3{\cdot}1$ |
| $2^5 3{\cdot}5^2$ | $-7^4$ | $1$ |
| $-113^2$ | $7^4$ | $2^7 3^4{\cdot}1$ |
| $-2^2{\cdot}61^2$ | $11^4$ | $3^5{\cdot}1$ |
| $239^2$ | $-2{\cdot}13^4$ | $1$ |
| $287^2$ | $-17^4$ | $2^7 3^2{\cdot}1$ |
| $2{\cdot}5861^2$ | $-3^2{\cdot}59^4$ | $7^9$ |
| $2{\cdot}13^2$ | $-3{\cdot}19^4$ | $5^8$ |
| $437147^2$ | $-21769^3$ | $2^9 3^4{\cdot}5^{12}$ |
| $1169^2$ | $2^3 3^4{\cdot}5^4$ | $-11^6$ |
| $2591^2$ | $-3{\cdot}43^4$ | $2{\cdot}11^6$ |
| $14089^2$ | $-131^4$ | $2^{11} 3{\cdot}5^6$ |
| $3^2{\cdot}5^3$ | $2^2{\cdot}19^3$ | $-13^4$ |
| $-37^3$ | $2{\cdot}29^3$ | $3{\cdot}5^4$ |
| $-2^8{\cdot}7^3$ | $3{\cdot}29^3$ | $11^4$ |
| $-71^3$ | $23^3$ | $2^4 3^2{\cdot}7^4$ |
| $-2{\cdot}203^3$ | $3^3{\cdot}79^3$ | $43^4$ |
| $2293^2$ | $2{\cdot}67^3$ | $-3{\cdot}5^9$ |
| $1079^2$ | $-2^{10}{\cdot}19^3$ | $3{\cdot}5^9$ |
| $-3{\cdot}36553^2$ | $203^3$ | $2^{11}{\cdot}5^9$ |
| $-138743^2$ | $3{\cdot}1027^3$ | $2^{13}{\cdot}5^9$ |
| $2^3 3^2{\cdot}12515^2$ | $-2797^3$ | $13^9$ |
| $107567^2$ | $-3{\cdot}3155^3$ | $2^{11}{\cdot}7^9$ |

| $a{\cdot}x^2$ | $b{\cdot}y^3$ | $c{\cdot}z^8$ |
|---|---|---|
| $-2 \cdot 3{\cdot}263^2$ | $29^3$ | $5^8$ |
| $-2^2{\cdot}353^2$ | $3^4{\cdot}11^3$ | $5^8$ |
| $5239^2$ | $-3^2{\cdot}163^3$ | $2{\cdot}7^8$ |
| $-2^5{\cdot}4015^2$ | $799^3$ | $7^8$ |
| $26311^2$ | $2^4 3^2{\cdot}95^3$ | $-13^8$ |
| $-32039^2$ | $241^3$ | $2^5 3^4{\cdot}5^8$ |
| $-39313^2$ | $2{\cdot}767^3$ | $3{\cdot}11^8$ |
| $31987^2$ | $-1489^3$ | $2^3 3^6{\cdot}5^8$ |
| $36631^2$ | $1679^3$ | $-2^6 3^5{\cdot}5^8$ |
| $-99431^2$ | $2^3{\cdot}1073^3$ | $3^2{\cdot}5^8$ |
| $-160975^2$ | $2^{10} 3{\cdot}203^3$ | $11^8$ |
| $-185039^2$ | $1633^3$ | $2^6 3^4{\cdot}7^8$ |
| $2^2{\cdot}774517^2$ | $-15613^3$ | $3^8{\cdot}11^8$ |
| $-6827035^2$ | $2{\cdot}28559^3$ | $3^3{\cdot}13^8$ |
| $9101359^2$ | $-43873^3$ | $2^7 3^7{\cdot}7^8$ |
| $-26615519^2$ | $78913^3$ | $2^7 3^5{\cdot}17^8$ |
| $-30042907^2$ | $2^3 3^3{\cdot}16037^3$ | $43^8$ |
| $2{\cdot}45707519^2$ | $1171537^3$ | $-3^5{\cdot}95^8$ |

These maps themselves are three point covers. In fact $f_3$ and $f_4$ are equations describing how the modular curve $X_0(N)$ covers the $j$-line $X_0(1)$ for $N = 2, 3$.

The base-change operation is quite efficient. Not using [Cog], but rather starting from just $\tau = -2$ corresponding to $2 - 3 + 1 = 0$, one gets 73 of the 101 orbits. The computer search through cutoff $10^9$ yields only 18 more orbits. Representatives of two of these appear on the left side of 3.4 and 3.5. Applying the base-change operations again yields nine more orbits. Note that the orbits appearing in 3.4 and the right side of 3.5 are three of the 10 known members of $T_H(\mathbf{Z})$; the orbit on the left side of 3.5 is particularly interesting too, being in $T_{3,3,4}(\mathbf{Z}^{\{2\}})$.

$$f_3(3^5,\ -2^2 61^2,\ 11^4) \;=\; (-13^3 1201^3,\ 2^2 61^2 12697^2,\ 3^8 11^8) \tag{3.4}$$

$$f_2(3^3 79^3,\ -2^1 7^3 29^3,\ 43^4) \;=\; (2^3 3^3 7^3 29^3 79^3,\ -109^2 275623^2,\ 43^8) \tag{3.5}$$

A computer search through cutoff $10^{11}$, complete for the easier cases $(2, \geq 3, \geq 8)$ found the remaining three orbits, one of which, $(-2 \cdot 184211^2,\ 4079^3,\ 3)$, is the only entry of Table 3.1 not picked up by the search previously.

Table 3.3. Order of the specialization sets $T^*_{p,q,r}$

| | | Section: | 4 | 5 | 6 | 7 | 8 | 9 |
|---|---|---|---|---|---|---|---|---|
| $p$ $q$ $r$ | | $\|T^*_{p,q,r}\|$ | $A_6$ | $SL_2(8)$ | $SL_3(3)$ | $PSp_4(3)$ | $A_9$ | $Sp_6(2)$ |
| 2 3 8 | | 99 | $f_{10}$ | | | | | |
| 2 3 9 | | 87 | | $f_9$ | | | | |
| 2 3 12 | | 82 | | | | | | |
| 2 3 $\geq 16$ | | 81 | | $f_{18}$ | | | | |
| 2 4 6 | | 48 | | | | | | |
| 2 4 8 | | 45 | | | $f_{26}$ | | | |
| 2 4 9 | | 45 | | | | $f_{27d}$ | | |
| 2 4 $\geq 12$ | | 44 | | | | | | |
| 2 6 | | 36 | | | | | | |
| 2 $\geq 8$ | | 35 | | | | $f_{27abc}$ | $f_{9,1}$ | $f_{28}$ |
| 3 3 4 | | 39 | $f_6$ | | | | | |
| 3 3 $\geq 6$ | | 27 | | | $f_{13b}$ | | | |
| 3 4 | | 24 | | | $f_{13c}$ | | | |
| 3 $\geq 6$ | | 23 | | | | | | |
| $\geq 4$ | | 21 | | | | | | |

In the rest of this paper we abbreviate the part of $T_{p,q,r}(\mathbf{Z}^{\{2,3\}})$ appearing on Tables 3.1 and 3.2 by $T^*_{p,q,r}$. Table 3.3 gives the cardinality of $T^*_{p,q,r}$ in all cases, under the normalization hypothesis $p \leq q \leq r$. The block on the right illustrates that we are using the bulk of Tables 3.1 and 3.2 in Sections 4-9.

## 4. $P\Gamma L_2(9)$ and index two subgroups; basics of computations

In this section we start from the following two covers.

$$\begin{aligned} \Lambda_6 &= (3A,\ 3B,\ 4A) \to (33,\ 3111,\ 42) \\ f_6(t,x) &= (x^2-2)^3 + t(3x-4)^2 \\ D_6(t) &= 2^{13}3^6t^4(t-1)^2 \end{aligned}$$

$$\begin{aligned} \Lambda_{10} &= (3AB,\ 2D,\ 8A) \to (3331, 22222, 811) \\ f_{10}(t,x) &= (x^3+12x^2+60x+96)^3x + 1728t(3x^2+28x+108) \\ D_{10}(t) &= -2^{99}3^{42}t^6(t-1)^5 \end{aligned}$$

The component set $C_6 = \pi_0(X_6(\mathbf{C}))^{\mathrm{cut}}$ of the cover $X_6 \to \mathbf{P}^1_{\mathbf{Q}}$ has six elements. The monodromy transformations $m_{6,0}$, $m_{6,1}$, and $m_{6,\infty}$ generate the alternating subgroup $\mathrm{Alt}(C_6)$. The Galois group is all of $\mathrm{Sym}(C_6)$, the constant field extension with Galois group $G_6/M_6$ is $\mathbf{Q}(\sqrt{2})$, as can be seen from considering $D_6(t)$ modulo squares. Similarly, $C_{10}$ can be identified with a projective line over a finite field $\mathbf{F}_9$ so that the monodromy group $Mo_{10}$ is $PGL_2(9)$ and the Galois group $G_{10}$ is $P\Gamma L_2(9)$; the constant field extension is again $\mathbf{Q}(\sqrt{2})$, although this can no longer be seen from $D_{10}(t)$. (In this case we write $Mo_{10}$ for the monodromy group, since $M_{10}$ is the standard notation for a different index-two subgroup of $P\Gamma L_2(9)$.)

The covers in the later section will be presented in the same three-line format, the meaning of the first line being as follows. This line gives first the conjugacy class $[m_c]$ of $m_c$ in $M$, in Atlas notation, for $c = 0, 1, \infty$. When $M$ has non-trivial outer automorphisms, like in the two cases here, there may be some ambiguity in how Atlas notation is used; here for example we could just as well write $8B$ rather than $8A$. The first line gives next the orbit-partition $\lambda_c$ of $m_c$ acting on the component set $C$. All the covers in this paper are rigid, meaning that they are completely determined by the group-theoretical data $([m_0]; [m_1], [m_\infty]; M \subseteq \mathrm{Sym}(C))$. The computation of the defining equation starts from this group-theoretical information.

Except for $f_{9a}$ in Section 5 and the $f^*_{N,m}$ in Section 8, all the defining equations have the simple form $f(t,x) = A(x) + tC(x)$. In these cases $x$ is a coordinate on the covering curve $X$, identifying it with $\mathbf{P}^1$. The cover $X \to \mathbf{P}^1$ is given by the rational function $x \mapsto -A(x)/C(x)$. Define $B(x)$ by $A(x) + B(x) + C(x) = 0$. Then the functions $A(x)$, $B(x)$, and $C(x)$ factor over $\mathbf{Q}$ according to $\lambda_0$, $\lambda_1$, and $\lambda_\infty$. Thus, for example, $f_{10}(1,x)$ factors as a quintic squared, according to $\lambda_{10,1} = 22222$. In brief, one determines the coefficients of $A$ and $C$ by imposing some normalization conditions and demanding that $B$ factor in the proper way. For very simple computations of this form see [Bir]; to make the computations feasible one has to use the "differentiation trick" described there systematically.

We specialize the cover $X_6$ at the 39-element set $T^*_{3,3,4}$. The 39 algebras $L_{6,\tau}$ are all of degree six, thus the factor fields all appear in [JR1]. The algebra $L_{6,\tau}$ is the twin [Rob1] of the algebra $L_{6,1-\tau}$, so that their splitting fields in $\mathbf{C}$ are identical: $K_{6,\tau} = K_{6,1-\tau}$. Since $L_{6,1/2}$ is self-twin, it is forced to be a non-field; in fact $L_{6,1/2}$ factors as a quintic field $L_5$ times $\mathbf{Q}$, with $L_5$ the unique quintic $\{2,3\}$-field with Galois group the Frobenius group $F_5$.

Here the 39 specialization points yield 37 isomorphism classes of algebras, the repetitions being $L_{6,2^2 19^3/13^4} \cong L_{6,2^2}$ and $L_{6,3^2 5^3/13^4} \cong L_{6,-3}$. Excluding $L_5 \times \mathbf{Q}$ and counting the repeated fields once, one has 10 $S_6$ twin pairs and 2 $C_3^2.D_4$ twin pairs. The remaining 6 twin pairs consist of a field and a non-field. Two of these have Galois group $S_5$, three have Galois group $S_4 \times S_2$, and one has Galois group $S_3 \times S_3$.

Note from Table 1.1 that the sets $NF(\{2,3\}, S_5)$ and especially $NF(\{2,3\}, S_6)$ are "unbalanced" with respect to the natural partition indexed by discriminant classes $d \in \{-6,-3,-2,-1,2,3,6\}$. In each case just under half of the fields have discriminant class 2. So one can view the existence of the cover $f_6$ as explaining this unbalance; it accounts for both $S_5$'s and 10 of the 13 $S_6$'s.

We specialize the cover $X_{10}$ at the 99-element set $T^*_{3,2,8}$. There is a tight relation between the covers $X_6$ and $X_{10}$. Namely $f_6(s,x)$, $f_6(1-s,x)$ and $f_{10}(1-(2s-1)^2, x)$ have the same splitting field over $\mathbf{Q}(s)$. This accounts for the behavior of the nineteen specialization points with $1 - \tau \in \mathbf{Q}^{\times 2}$ .

Twelve of the elements $\tau \in T^*_{3,2,8}$ satisfy $2(1-\tau) \in \mathbf{Q}^{\times 2}$. By the form of $D_{10}(t)$, the corresponding $G_\tau$ are in the index two subgroup $M_{10}$ of $P\Gamma L_2(9)$. The polynomials $f_{10}(-23^3, x)$ and $f_{10}(-17^3 47^3/7^8, x)$ each factor as $f_9 f_1$. The field $K_{10,-23^3}$ is one member of the two-element set $NF(\{2,3\}, C_3^2.C_4)$. The field $K_{10,-17^3 47^3/7^8}$ has Galois group $C_3^2.Q$, $Q$ being the quaternion group. One knows that $NF(\{2,3\}, Q)$ also has two elements, with defining equations

$$g_\pm(x) = x^8 \pm 12x^6 + 36x^4 \pm 36x^2 + 9.$$

The field $K_-$ is totally real while $K_+$ is totally imaginary; the quaternionic subfield of $K_{10,-17^3 47^3/7^8}$ is $K_+$. The remaining ten fields are all distinct, with Galois group $M_{10}$, and discriminants $2^a 3^b$, $24 \le a \le 34$, $10 \le b \le 18$.

For eight of the remaining sixty-eight $\tau$, the polynomial $f_{10}(\tau, x)$ factors as $f_9 f_1$, namely $\tau = -2^1 61^3$, $-2^1 28559^3/3^3 13^8$, $-239^3/3^{13}$, $-213^3 59^3/3^1 11^8$, $-29^3/5^8$, $2/3^3$, $311^3 3767^3/3^5 5^8 19^8$, and $-4079^3/3$. The splitting fields $K_{10,\tau}$ all have Galois group of the form $C_3^2.B$, with $B$ the Sylow 2 subgroup of $GL_2(3)$; these eight fields are all distinct. The remaining sixty specialization points give 59 fields in $NF(\{2,3\}, P\Gamma L_2(9))$, the unique duplication being $K_{10,13^3 1201^3/3^8 11^8} = K_{10,5^3/2^2}$, verified by *polred*. Note that none of the 59 algebras $L_{10,\tau}$ are tame at 2, because all of them have $\mathbf{Q}(\sqrt{2})$ in their splitting field. However exactly one is tame at 3, namely

$$L_{10,11^3/2^3} : x^{10} - 4x^9 + 6x^8 + 24x^2 + 32x + 16, \tag{4.1}$$

with discriminant $-2^{34}3^9$. Note that $f(11^3/2^3, x)$ has relatively large coefficients. The displayed polynomial is obtained from $f(11^3/2^3, x)$ by *polred*. We have used *polred* analogously below, so as to always display only polynomials with relatively small coefficients.

We now discuss Frobenius computations somewhat, as promised in Section 2. Let $G^\natural$ be the set of conjugacy classes in $G = P\Gamma L_2(9)$. For every $p \ge 5$ one has a Frobenius element $\mathrm{Fr}_{\tau,p}$ in $G^\natural$. As an element of the quotient group $C_2 \times C_2$, it is given by a pair of given by a pair of Jacobi symbols:

$$[\mathrm{Fr}_{\tau,p}] = \left( \left( \frac{2}{p} \right), \left( \frac{2(1-\tau)}{p} \right) \right).$$

The Frobenius class $\mathrm{Fr}_{\tau,p}$ itself is then completely determined by the factorization pattern of $f_{10}(\tau, x)$ in $\mathbf{Q}_p$. The possibilities are

| | $A_6 = PSL_2(\mathbf{F}_9) = M_{10}'$ | | | | | $S_6 - A_6$ | | | $PGL_2(9) - PSL_2(9)$ | | | $M_{10} - M_{10}'$ | |
|---|---|---|---|---|---|---|---|---|---|---|---|---|---|
| $[\mathrm{Fr}_{\tau,p}]$ | $(1,1)$ | | | | | $(-1,-1)$ | | | $(1,-1)$ | | | $(-1,1)$ | |
| Atlas | $1A$ | $2A$ | $3AB$ | $4A$ | $5AB$ | $2BC$ | $4B$ | $6AB$ | $2D$ | $8AB$ | $10AB$ | $4C$ | $8CD$ |
| # | 1 | 45 | 80 | 90 | 144 | 30 | 90 | 240 | 36 | 180 | 144 | 180 | 180 |
| $\mathrm{Fr}_{\tau,p}$ | $1^{10}$ | $2^4 1^2$ | $3^3 1$ | $4^2 1^2$ | $5^2$ | $2^3 1^4$ | $4^2 2$ | $631$ | $2^5$ | $81^2$ | $10$ | $4^2 1^2$ | $82$ |

Consider the sixty $\tau$ mentioned last above. The groups $G_\tau$ surject onto $C_2 \times C_2$, exactly because $2(1-\tau) \ne 1, 2 \in \mathbf{Q}^{\times 2}$. From the list of maximal subgroups of $P\Gamma L_2(9)$ in the Atlas, one knows that any proper such group is divisible by at most two of $\{2,3,5\}$. Working with just the 20 primes $p = 5, 7, \ldots, 79$, one finds that 30 divides $|G_\tau|$ always, and so $G_\tau = P\Gamma L_2(9)$ always. Even the first ten primes suffice to prove that the fields are all pairwise non-isomorphic, except for perhaps the exception noted above.

In practice we usually work with larger primes, so that Frobenius elements can be computed by factoring merely over $\mathbf{F}_p$. Also knowing the generic relative frequencies, as indicated by # above, is an aid in deciding when Frobenius data is truly suggestive of nongeneric behavior. In practice, we work with a large number of primes, so that Frobenius data can be expected to very clearly distinguish between generic and non-generic behavior.

## 5. $SL_2(8)$ and $\Sigma L_2(8)$; Shimura curves

Here we specialize the two covers. This first is

$$\begin{aligned}
\Lambda_{9a} &= (333, 22221, 9) \\
f_{9a}(t,x) &= (x^3 - 9x^2 - 69x - 123)^3 - \\
&\quad 2^{14}t(9x^4 - 42x^3 - 675x^2 - 1485x - 441) - 2^{28}t^2 \\
D_{9a}(t) &= 2^{140}3^{18}t^6(t-1)^4 \cdot \\
&\quad (4398046511104t^3 - 5421322469376t^2 - 7496810496t + 513922401)^2.
\end{aligned}$$

The defining polynomial $f_{9a}$ was sent to us by Elkies in 1995. This cover has genus one, while all the other covers in this paper have genus zero. This is the reason behind the extraneous cubic-squared factor in $D_{9a}(t)$.

To get the second cover, put

$$\begin{aligned}
A(x) &= x^9 + 108x^7 + 216x^6 + 4374x^5 + \\
&\quad 13608x^4 + 99468x^3 + 215784x^2 + 998001x + 810648.
\end{aligned}$$

In [Mat, page 193], Matzat introduces the cover

$$\begin{aligned}
\Lambda_{9b} &= (3B, (3B)^2, 9A) \to (33111, 33111, 9) \\
f_{9b}(u,x) &= A(x) + u2^{13}3^4 \\
D_{9b}(u) &= 2^{104}3^{50}(u^2+3)^4.
\end{aligned}$$

Its monodromy and Galois group are both $\Sigma L_2(8)$. The resolvent $A_3$ cubic is itself a cover of the $u$-line being

$$\begin{aligned}
\Lambda_3 &= (3A, (3A)^2, 1A) \to (3, 3, 111) \\
f_3(u,x) &= x^3 - (9u^2+27)x - (9u^3 + 9u^2 + 27u + 27) \\
D_3(u) &= 3^6(u^2+3)^2(u-3)^2.
\end{aligned}$$

Matzat's three-point cover does not fit into our set-up because the ramification locus is $\{\sqrt{-3}, -\sqrt{-3}, \infty\}$. Putting $t = -u^2/3$, we work with the "double":

$$\begin{aligned}
\Lambda_{18} &= (2^9,\ 3^4\,1^6,\ 18) \\
f_{18}(t,x) &= f_{9b}(u,x)f_{9b}(-u,x) \\
&= A(x)^2 + 2^{26}3^9t \\
D_{18}(t) &= -2^{460}3^{189}t^9(t-1)^8.
\end{aligned}$$

The Galois group is $G_{18} = \Sigma L_2(8)^2.2$. The monodromy group is the unique index three normal subgroup. The constant field extension, with Galois group $G_{18}/M_{18}$, is $\mathbf{Q}(2\cos(2\pi/9))$, with defining equation $x^3 - 3x + 1$; this is the unique field in $NF(\{2,3\}, A_3)$.

We specialize $f_{9a}$ at the 87-element set $T^*_{3,2,9}$. The polynomial $f_{9a}(-17^3/2^7, x)$ factors as $f_8 f_1$, the degree eight factor having Galois group $A_4 \times C_2$. The polynomials $f_{9a}(-2^5 3^2, x)$ and $f_{9a}(3^2 17^3/2^{15}, x)$ define the same field, with defining equation $x^9 + 6x^3 - 2$ and Galois group the 54-element affine group $(\mathbf{Z}/9).(\mathbf{Z}/9)^\times$. The remaining 84 specialization points give 55 fields, all with Galois group $\Sigma L_2(8)$. All these fields are wildly ramified at three. Exactly two of them are only tamely ramified at 2:

$$L_{9a,\tau} \ : \ x^9 - 36x^6 - 162x^4 - 54x^3 - 972x^2 + 486x - 594 \tag{5.1}$$

$$L_{9a,4/3} \ = \ x^9 - 18x^3 + 27x - 6. \tag{5.2}$$

The first field arises from $\tau \in \{2^5/3^4, 2^{10}19^3/5^93, -2^95^3/3^2, -73^3/2^{15}3^2\}$. The field discriminants are $2^83^{26}$ and $2^63^{26}$, respectively.

We specialize $f_{18}$ at the 81-element set $T^*_{2,3,18} = T^*_{2,3,\infty}$. When $-3\tau$ is not a square in $\mathbf{Q}$ we get fields which generically have Galois group containing two copies of $SL_2(8)$, thus interesting, but out of our self-imposed context. The six elements $\tau$ with $-3\tau$ a square are $\tau = -u^2/3$ with $u \in \{1, 10/9, 35/18, 3, 595/108, 37\}$. The twelve algebras $L_{9b,\pm u}$ are all fields, and so is $L_{9b,0}$; they are all pairwise non-isomorphic. Exactly one of these fields appeared already as a specialization of $f_{9a}$, namely $L_{9b,-5\cdot 7/2\cdot 3^2}$ which coincides with the tame-at-two field $L_{9a,4/3}$. Of the twelve remaining fields, nine have Galois group $\Sigma L_2(8)$, and three have Galois group $SL_2(8)$, these being

$$\begin{array}{rcl} L_{9b,-3} & : & x^9 - 12x^6 - 18x^5 + 36x^2 - 27x - 128 \\ L_{9b,37} & : & x^9 - 36x^6 - 54x^5 - 324x^4 - 216x^3 - 972x^2 - 243x - 2124 \\ L_{9b,1} & : & x^9 - 36x^6 - 54x^5 + 432x^3 + 324x^2 - 243x - 1152. \end{array}$$

These fields have discriminants $2^{14}3^{22}$, $2^{14}3^{26}$ and $2^{14}3^{26}$ respectively.

In comparison with the other examples in this paper, the most remarkable phenomenon here is that often two or more $L_{9a,\tau}$ are isomorphic. This phenomenon can be partially explained as follows. The cuspidal data identifies $\mathbf{P}^1$ with a minimal-area Shimura curve $X_0(1)$ associated to the cubic field $\mathbf{Q}(\cos(2\pi/9))$ and no ramification at finite places. The cover $X_{9b} \to \mathbf{P}^1$ is identified with $X_0(2) \to X_0(1)$, the ideal $(2)$ having residual cardinality 8. But now also one has a degree four cover $\pi : X_0(P) \to X_0(1)$, where $P$ is the unique prime above 3. The curve $X_0(P)$ has an Atkin-Lehner-type involution $W_P$, and hence a second natural map $\pi \circ W_P$ to the base curve $X_0(1)$. As in the classical case, one can think of $X_0(P) \subset X_0(1) \times X_0(1)$ as being a correspondence $T_P$ from $X_0(1)$ to $X_0(1)$ of bidegree $(4, 4)$. In terms of our fixed coordinate $t$ on $X_0(1)$, and the same coordinate $s$ on the second copy of $X_0(1)$, the defining equation for $X_0(P)$ is unique up to scalars:

$$\begin{aligned} h_P(s,t) = \\ & 2^{12}st\,\big(-3^{10}17^366383 + 2^{13}3^81054805(s+t) + 2^{27}3^5211(s^2+t^2) + \\ & 2^{13}3^62486119st - 2^{21}3^322267(s^2t + st^2) - 2^{24}55s^2t^2 + 2^{30}(s^3t^2 + s^2t^3)\big) \\ & -3^6(-3^217^3 + 2^{15}s)^3 - 3^6(-3^217^3 + 2^{15}t)^3 - 3^{12}17^9. \end{aligned}$$

The specialization points giving the same field are exactly as follows. Here $\tau \approx \sigma$ rather than just $\tau \sim \sigma$, indicates that $h_P(\sigma, \tau) = 0$; so these isomorphisms $L_{9a,\sigma} \cong L_{9a,\tau}$ are explained by the Hecke correspondence $T_P$.

$$\begin{array}{rcl} -1 & \sim & 2797^3/13^9 \\ 5^37^3/3^5 & \sim & -4079^3/3 \\ 3^2 & \sim & 5^3/3^3 \\ 3^5 & \sim & -23^3 \\ -2^53^2 & \sim & 3^217^3/2^{15} \\ -1/2^3 & \sim & 3^567^3/2^3 \end{array} \qquad \begin{array}{rcccc} -2^4/3^2 & \sim & 2^27^3/3 & \approx & 23^3/2^{12}3 \\ 2^3/3^2 & \sim & -13^3/2^23 & \approx & -5^6/2^93 \\ -2^33 & \sim & 1/2^2 & \approx & 3^111^3/2^9 \\ -3 & \sim & 5^3/2^7 & \approx & 3^119^3/2^7 \end{array}$$

$$\begin{array}{rcccccccc} & & 2^5/3^4 & \sim & 2^{10}19^3/5^93 & \sim & -2^95^3/3^2 & \approx & -73^3/2^{15}3^2 \\ -2^3 & \sim & -3^2/2^4 & \sim & 5^3/2^2 & \approx & 7^3/2^9 & \approx & 3^273^3/2^9 \\ 3/2^2 & \sim & -11^3/2^8 & \approx & -5^37^33/2^8 & \sim & -7^329^3/2^{11}5^9 & \approx & 5^3631^33/1^{11}7^9 \\ -2^43 & \sim & 2^2 & \approx & 3^15^3/2^{12} & \sim & -5^3383^3/2^{13} & \approx & -3^113^379^3/2^{13}5^9 \end{array}$$

The Hecke correspondence $T_P$ also explains two of the three cases mentioned above of unexpectedly small Galois group. First, the degree seven polynomial $h_P(t,t)$ factors as $f_3^2 f_1$, the root of the linear factor being $\tau = -17^3/2^7$. So this $\tau$ is a CM point, which forces $L_{9a,\tau}$ to be non-generic. Second, the degree four polynomial $h_P(0,t)$ factors as $f_3 f_1$, the root of the linear factor being $\tau = 3^2 17^3/2^{15}$. So this $\tau$, like 0, is also a CM point, again forcing non-genericity in $L_{9a,\tau}$.

The cover $X_{10}$ from the previous section is also a minimal area Shimura curve, coming from the field $\mathbf{Q}(\sqrt{2})$. For more on Shimura curves as covers of the projective line, see [Tak], [Elk2].

## 6. $SL_3(3)$ and $SL_3(3).2$; projective twinning

Here we use two new covers $13b$, $13c$, and a cover 26, which is a doubled version of one of the covers in [Mal1]. The new covers are

$$\begin{aligned}
\Lambda_{13b} &= (3B,\ 3A,\ 8A) \to (33331,\ 3331111,\ 841) \\
f_{13b}(t,x) &= \left((x^3 - 6x + 6\sqrt{-2}x - 4 - 8\sqrt{-2})(x - 2 - \sqrt{-2})\right)^3 (x - 2 + 2\sqrt{-2}) - \\
&\quad t2^2 3^2 (3x - 4 + \sqrt{-2})^4 (3x - 8 + \sqrt{-2}) \\
D_{13c}(t) &= 2^{92} 3^{72} (1 + \sqrt{-2})^{96} t^8 (t-1)^6
\end{aligned}$$

$$\begin{aligned}
\Lambda_{13c} &= (4A,\ 3A,\ 8A) \to (44221,\ 3331111,\ 841) \\
f_{13c}(t,x) &= (x^2 - 3\sqrt{-2}x - 3\sqrt{-2} - 3)^4 (x^2 + 6\sqrt{-2} - 3)^2 (x + 3 - 3\sqrt{-2}) + \\
&\quad t2^2 3^3 (1 + \sqrt{-2})^8 (x + 1 - \sqrt{-2})^4 (3x + 5 - \sqrt{-2}) \\
D_{13c}(t) &= 2^{92} 3^{72} (1 + \sqrt{-2})^{96} t^8 (t-1)^6.
\end{aligned}$$

These covers were first computed over $\mathbf{F}_{11}$, where each cover appeared with its projective twin, as described below; at this level there is no evident relation between the defining equation for a cover ($\sqrt{-2} \mapsto 3$) and the defining equation for its twin ($\sqrt{-2} \mapsto 8$). We worked with the auxiliary prime 11, because it is the smallest prime besides 3 which is split in the field $\mathbf{Q}(\sqrt{-2})$. We then lifted via the 11-adics to solutions in $\mathbf{Q}(\sqrt{-2})$ as explained in [Mal2]. At this level passing from a cover to its twin is induced by complex conjugation in the ground field $\mathbf{Q}(\sqrt{-2})$.

Put

$$\begin{aligned}
A(x) &= (x^3 - 12x^2 - 6x - 64)(x^4 + 16x^3 - 36x^2 + 128x - 28) \cdot \\
&\quad (x^6 + 12x^5 + 54x^4 + 176x^3 + 444x^2 + 624x + 552) \\
B(x) &= (x^4 + 16x^3 + 72x^2 + 128x + 188)(3x^4 - 4x^3 + 12x^2 - 24x - 68)^2.
\end{aligned}$$

The cover appearing in [Mal1] is ramified at $(-\sqrt{-8}, \sqrt{-8}, \infty)$:

$$\begin{aligned}
\Lambda_{13a} &= (8A, 8B, 2A) \to (841, 841, 222211111) \\
f_{13a}(u,x) &= A(x) - uB(x) \\
D_{13a}(u) &= 2^{160} 3^{114} (u^2 + 8)^{10}.
\end{aligned}$$

One has $M_{13a} = G_{13a} = SL_3(3)$. Just as we doubled the Matzat cover, we now put $t = -u^2/8$, and double the Malle cover to put it into our context:

$$\begin{aligned}
\Lambda_{26a} &= (2B,\ 8AB,\ 4B) \to (2^{13},\ 8^2\,4^2\,1^2,\ 4^4\,2^5) \\
f_{26a}(t,x) &= f_{13a}(u,x) f_{13a}(-u,x) \\
&= A(x)^2 + 8tB(x)^2 \\
D_{26a}(t) &= -2^{777}3^{452}t^{13}(t-1)^{20}.
\end{aligned}$$

In distinction to what happened when we doubled the Matzat cover, here the groups have not gotten much bigger: $M_{26} = G_{26} = SL_3(3).2$.

To facilitate comparison with $f_{26a}$ below, we consider the degree 26 polynomials $f_{26b} = f_{13b}\overline{f}_{13b}$ and $f_{26c} = f_{13c}\overline{f}_{13c}$ in $\mathbf{Q}[t,x]$. We specialize $f_{26b}$ at the 27-element set $T^*_{3,3,8}$ and $f_{26c}$ at the 24-element set $T_{4,3,8}$. The polynomials $f_{26b}(-8,x)$, $f_{26c}(4,x)$, $f_{26c}(3/4,x)$ each factor, and the corresponding splitting fields have Galois group $3^{1+2}_+.D_4$, $3^{1+2}_+.V$, and $3^{1+2}_+.C_2$ respectively. A Frobenius computation shows that the remaining 48 polynomials have Galois group all of $SL_3(3).2$. The Frobenius computation shows that these fields are all non-isomorphic except for perhaps $L_{26c,1/4}$ and $L_{26c,-8}$. We have verified that this last pair of fields is indeed isomorphic. We used the method described in Section 7, working only with the roots of $f_{13c,1/4}$ and $f_{13c,-8}$.

Now consider cover 26a. In general, if $-2\tau$ is not a square in $\mathbf{Q}$ then $L_{26a,\tau}$ contains a subfield, isomorphic to $\mathbf{Q}(\sqrt{-2\tau})$. If $-2\tau$ is a square then one has a factorization $f_{26a}(\tau,x) = f_{13a}(u,x)f_{13a}(-u,x)$, wth $u^2 = -8\tau$. Thirty-eight of the elements in the 45-element set $T^*_{2,8,4}$ are such that $-2\tau$ is a non-square. A Frobenius computation shows that all of these have Galois group all of $SL_3(3).2$ and are pairwise non-isomorphic. Moreover the fields $L_{26a,\tau}$ are not isomorphic with any of the 47 fields coming from $26b$ and $26c$.

The seven elements $\tau$ with $-2\tau$ a square are $\tau = -u^2/8$ with $u = 1, 2, 7/2, 4, 8$, 44 and 10. A Frobenius computation shows generic behavior in the first six cases: the splitting field of $f_{13a}(u,x)$ has Galois group all of $SL_3(3)$, and these six splitting fields are non-isomorphic. Note that the degree thirteen fields $\mathbf{Q}[x]/f_{13a}(u,x)$ and $\mathbf{Q}[x]/f_{13a}(-u,x)$ are non-isomorphic, corresponding to different permutation representations of their common Galois group. This is seen clearly in the degenerate case $\tau = 5^2/2$. In this case $f_{13a}(-10,x)$ factors as $9+4$ and $f_{13a}(10,x)$ factors as $12+1$; the common Galois group is $3^{1+2}_+.\tilde{S}_4$.

The twinning phenomenon here is quite general, deriving from the fact that the permutation representations of a group $PGL_n$ on the projective space $\mathbf{P}^{n-1}$ and its dual $\check{\mathbf{P}}^{n-1}$ are not isomorphic, for $n \geq 3$. A prime $\ell$ and a degree $n \geq 3$ being fixed, say that a number field of degree $\ell^{n-1} + \ell^{n-2} + \cdots + \ell + 1$ is projective iff its Galois group can be identified with $PGL_n(\ell)$, as a permutation group. Then a projective field $L$ has a non-isomorphic twin $L^t$, there being a canonical isomorphism $L = L^{tt}$.

This situation of twin projective fields contrasts in one important way with the sitation of twin $A_6$ or $S_6$ sextic fields, which played a role in Section 4. For, only in the projective case, the two different permutation representations induce the same linear representation. So, only in the projective case, do twin fields have the same Dedekind zeta function. One can summarize by saying that $A_6$ and $S_6$ sextic fields come in fraternal twin pairs while projective fields come in identical twin pairs. In particular, twin projective fields have the same discriminant. In the cases here, one

has discriminants $2^{28}3^{24}$, $2^{36}3^{26}$, $2^{28}3^{22}$, $2^{38}3^{26}$, $2^{32}3^{24}$, $2^{38}3^{20}$ for $|s| = 1, 2, 7/2, 4, 8, 44$.

## 7. $PSp_4(3)$ and $SO_5(3)$; exhibiting exceptional isomorphisms

Here we work with four degree twenty-seven covers, denoted $27a$, $27b$, $27c$, $27d$. Cover 27a is due to Häfner [Häf], while 27b and 27c are new here. These three covers are related to algebraic hypergeometric functions with finite monodromy, corresponding to the entries 47, 45, and 48 on the Beukers-Heckmann list [BH]; they all have $M = G = SO_5(3) = W(E_6)$. Cover 27d is new also; it is related to hypergeometric-like functions with finite monodromy, as explained in [Rob2]. For this cover, $M = PSp_4(3)$ and $G = PGp_4(3) = SO_5(3)$, the corresponding constant field extension being $\mathbf{Q}(\sqrt{-3})$.

$$\begin{aligned}
\Lambda_{27a} &= (9AB,\ 2C,\ 12C) \to (9^3,\ 2^6\,1^{15},\ 12\,6\,4^2\,1)\\
f_{27a}(t,x) &= (x^3+6\,x^2-8)^9 - t2^4 3^{12} x^6 (x+5\,x+4)^4(x-2)\\
D_{27a}(t) &= 2^{414}3^{450}t^{24}(t-1)^6
\end{aligned}$$

$$\begin{aligned}
\Lambda_{27b} &= (12AB,\ 2C,\ 8A) \to (12^2\,3, 2^6\,1^{15},\ 8^3\,2\,1)\\
f_{27b}(t,x) &= 2^4x^3(x^2-3)^{12} - t3^9(x-1)(x-2)^8(x^2-2x-1)^8\\
D_{27b}(t) &= 2^{542}3^{270}t^{24}(t-1)^6
\end{aligned}$$

$$\begin{aligned}
\Lambda_{27c} &= (9AB,\ 2C,\ 8A) \to (9^3,\ 2^6\,1^{15},\ 8^3\,2\,1)\\
f_{27c}(t,x) &= 2^{18}(x^3+9\,x^2+6\,x+1)^9 - 3^{15}\,t\,x\,(2\,x+1)^8(x^2-2\,x-1)^8\\
D_{27c}(t) &= 2^{522}3^{450}t^{24}(t-1)^6
\end{aligned}$$

$$\begin{aligned}
\Lambda_{27d} &= (9A,\ 2B,\ 4A) \to (9^3,\ 2^{10}\,1^7,\ 4^6\,1^3)\\
f_{27d}(t,x) &= 4(3x^3-12\,x-8)^9 - t(12\,x^3+54\,x^2+62\,x+22)\cdot\\
&\quad (9\,x^6+81\,x^5+135\,x^4+276\,x^3+432x^2+288\,x+64)^6\\
D_{27d}(t) &= 2^{520}3^{459}t^{24}(t-1)^{11}
\end{aligned}$$

The three new covers were all found first modulo 5, and then lifted via the 5-adics to rational solutions as explained in [Mal2]. The cover $X_{28}$ of Section 9 was also computed in this way. Of these four new covers, 27d was the hardest to obtain, since the $N+2$ parts of $3N$ are more evenly distributed among the cusps.

In cases $27a$, $27b$, and $27c$ we specialize to $\tau$ in the 35-element set $T^*_{9,2,12} = T^*_{12,2,8} = T^*_{9,2,8} = T^*_{\infty,2,\infty}$. In case $27d$ we specialize to $\tau$ in the 45-element set $T^*_{9,2,4}$. All 150 defining polynomials $f_{27*}(\tau,x)$ are irreducible, so all the algebras $L_{27*,\tau}$ are fields. A Frobenius computation proves that at least 146 of these fields have Galois group all of $G$ or $G'$. Using the monodromy technique described in Section 9, we have verified that each apparent group-drop indeed occurs, and computed degree nine resolvents as follows.

$$h_{27a,-1}(x) = x^9 - 3x^8 + 6x^7 - 3x^5 + 15x^4 + 6x^3 + 6x^2 + 6x + 2 \tag{7.1}$$

$$h_{27b,-1/48}(x) = x^9 - 12x^6 - 9x^5 + 54x^3 + 36x^2 + 18x - 56 \tag{7.2}$$

$$g_{27d,1/2}(x) = x^9 - 6x^7 + 27x^5 - 12x^4 - 174x^3 - 108x^2 + 183x + 124 \tag{7.3}$$

$$g_{27d,1/2\cdot 13^4}(x) = x^9 + 3x^8 - 16x^6 - 36x^5 + 36x^4 + 96x^3 - 108x - 36. \tag{7.4}$$

Here $h_{27a,-1}$ and $f_{27a,-1}$ have the same splitting field, and so too do $h_{27b,-1/48}$ and $f_{27b,-1/48}$; the Galois groups here have the form $3^3.(S_4 \times C_2)$ and $3^3.S_4$ respectively. On the other hand, the Galois groups $G_{27d,1/2}$ and $G_{27d,1/2\cdot 13^4}$ each have the form $3^{1+2}.\tilde{S}_4$; the corresponding degree nine resolvents have Galois group of the form $3^2.\tilde{S}_4$. The field discriminants of the four displayed nonic polynomials are $2^{13}3^{15}$, $2^{12}3^{18}$, $-2^{22}3^{13}$, and $-2^{22}3^{15}$ respectively.

The Frobenius computation proves that except for the possible isomorphism

$$L_{27b,-48} \cong L_{27d,32/81},$$

the 146 generic specializations are all distinct. To prove that $L_{27b,-48}$ and $L_{27d,32/81}$ are isomorphic one needs to exhibit the unique isomorphism; this is our topic for the rest of the section.

In numerical terms, one needs to find the unique bijection $\sigma$ from the complex roots $A = \{\alpha_i\}$ of $f_{27b,-48}$ to the complex roots $B = \{\beta_j\}$ of $f_{27d,32/81}$ such that the the polynomial

$$\text{test}_{\sigma,r}(x) = \prod_i (x - \alpha_i - r\beta_{\sigma(i)}) \tag{7.5}$$

has rational coefficients for all $r \in \mathbf{Q}$; here we introduce $r$ just to avoid possible inseparability problems; in practice one can simply fix $r$ at 1.

The bijection $\sigma$ must intertwine complex conjugation: $c_B \circ \sigma = \sigma \circ c_A$. Since the polynomials in question have exactly three real roots, there are $3!12!2^{12} \approx 1.177 \times 10^{13}$ such bijections. The sheer size of this number is why *Pari*'s *isisom*, which is designed to find isomorphisms between fields, does not work here. Indeed, if we were simply given the bare polynomials $f_{27b,-48}$ and $f_{27d,32/81}$, we would not know how to find the desired isomorphism.

But we are not given the polynomials in isolation. Rather we can figure out how the monodromy operators $m_c$ act on the roots. In the case of $f_{27b,-48}$ this action is as follows.

$$\begin{aligned} m_{27b,0} &= (1a, 15a, 7a, 4a, 3a, 2a, 2b, 3b, 4b, 7b, 15b, 1b)(5a, 6, 5b) \cdot \\ &\quad (11, 12b, 13b, 14b, 8b, 9b, 10, 9a, 8a, 14a, 13a, 12a) \\ m_{27b,\infty} &= (1a, 1b)(2a, 3a, 4a, 5a, 5b, 4b, 3b, 2b)(6, 7a, 8a, 9a, 10, 9b, 8b, 7b) \cdot \\ &\quad (11)(12a, 13a, 14a, 15a, 15b, 14b, 13b, 12b). \end{aligned}$$

Here we have ordered the $\text{Gal}(\mathbf{C}/\mathbf{R})$-orbits from left to right, 1 through 15. If a root is in the lower half plane, we append "a"; if it is the upper half plane we append "b"; it it is on the real line we index the root simply by a number; thus $\alpha_{1a} = -(12.976\ldots) - (243.333\ldots)i$ and $\alpha_6 = .518\ldots$. The action of $m_{27b,0}$ is obtained by considering the inverse image of $[-48, 0]$ in the complex plane with coordinate $x$. It is the union of three "wheels", two with 12 spokes and one with 3. The action is obtained by "rotating the wheels one spoke counter-clockwise". All this can be done numerically, working visually with say the roots of $f_{27b}(-48u^{12}, x)$, with $u = 0, 1/100, 2/100, \ldots, 1$; restricting to $u \geq 20/100$ lets one work entirely in standard precision and still do the computation. The other actions are obtained similarly, always in the the spirit of Grothendieck's *dessins d'enfants* [Sch]. The action of $m_{27b,\infty}$ is obtained similarly. Repeating the procedure for the roots of $f_{27d,32/81}$ gives

$$\begin{aligned} m_{27d,0} &= (1a, 2a, 10a, 6a, 3, 6b, 10b, 2b, 1b)(9b, 8b, 7b, 5b, 4, 5a, 7a, 8a, 9a)\cdot \\ &\quad (11a, 13a, 14a, 12a, 15, 12b, 14b, 13b, 11b) \\ m_{27d,1} &= (5a)(5b)(13a)(13b)(14a)(14b)(15)(1a, 1b)(2a, 12a)(2b, 12b)\cdot \\ &\quad (3, 4)(6a, 9a)(6b, 9b)(7a, 8a)(7b, 8b)(10a, 11a)(10b, 11b). \end{aligned}$$

Let $C$ be one of the two root sets $A$ or $B$. We use the monodromy action to view the 27-element set $C$ as a structured set. Namely consider the set $\mathrm{Sub}_2(C)$ of two-element subsets. Under the action of the monodromy group one has a decomposition into two orbits

$$\mathrm{Sub}_2(C) \quad = \quad \mathrm{Sub}_2(C)' \coprod \mathrm{Sub}_2(C)''.$$

Here the first orbit has 135 elements and the second has 216. In other words, we can view $C$ as the vertices of a graph $\Gamma_C$, with edge-set $\mathrm{Sub}_2(C)'$. The action of complex conjugation on $C$ extends to $\Gamma_C$. The graph-with-involution $(\Gamma_C, i_C)$ has only $2^7 3^2 = 1152$ automorphisms.

The bijection $\sigma : A \to B$ we seek respects not only complex conjugation but also the graph structure. To find it, we work purely group-theoretically to first find some isomorphism $s : (A, i_A) \to (B, i_B)$ respecting the graph structure; this is easy. Then we compose it with the 1152 elements of $\mathrm{Aut}(\Gamma_B, i_B)$ to get 1152 candidates for $\sigma$. For only one of them does 7.5 appear to have rational coefficients, and we have thereby numerically determined $\sigma$. It is as follows:

| roots of $f_{27b,-48}$: | $1a$ | $2a$ | $3a$ | $4a$ | $5a$ | $6$ | $7a$ | $8a$ | $9a$ | $10$ | $11$ | $12a$ | $13a$ | $14a$ | $15a$ |
|---|---|---|---|---|---|---|---|---|---|---|---|---|---|---|---|
| roots of $f_{27d,32/81}$: | $7b$ | $1b$ | $9a$ | $5b$ | $8a$ | $4$ | $6a$ | $2b$ | $11a$ | $15$ | $3$ | $14b$ | $13a$ | $12b$ | $10a$ |

Here since, e.g. $\sigma(\alpha_{1a}) = \beta_{7b}$ we must have $\sigma(\alpha_{1b}) = \beta_{7a}$ as well.

To check rigorously that indeed $\sigma$ induces an isomorphism we proceed as follows. First, we numerically solve the 27-by-27 system $Tz = B$ with $T_{ik} = \alpha_i^k$ and $B_i = \beta_{\sigma(i)}$ for the vector $z$, rationalizing at the end. Then

$$y = \sum_{k=0}^{26} z_k x^k, \tag{7.6}$$

as an element of the ring $\mathbf{Q}[x]/f_{27b}(-48, x)$, should satisfy $f_{27d}(32/81, y) = 0$. In our case 7.6 takes the explicit form $126414618624y =$

$$\begin{array}{rrr} -778448003x^{26} & -981509289x^{25} & -45939794218464x^{24} \\ 1137207587245554x^{23} & -12724591174373616x^{22} & 84720990963862440x^{21} \\ -370643301489686778x^{20} & 1104517184207752350x^{19} & -2221386267735948267x^{18} \\ 2770729912599438087x^{17} & -1364804638977353670x^{16} & -1659195683617968252x^{15} \\ 3294135546577106040x^{14} & -1424334213106643304x^{13} & -1494872066602993428x^{12} \\ 1744593287940021708x^{11} & 585147838091136 39x^{10} & -818860129717879323x^{9} \\ 171625804007741892x^{8} & 237402868177405098x^{7} & -685482891344617 68x^{6} \\ -50056832848713984x^{5} & 11561622988644846x^{4} & 7711324091701110x^{3} \\ -484480924209657x^{2} & -608102864271747x & -73397655884286. \end{array}$$

This indeed satisfies $f_{27d}(32/81, y) = 0$. The case $L_{13c,1/4} \cong L_{13c,-8}$ from the previous section is substantially easier, as the root sets to be identified have only 13 elements each.

## 8. $A_9$ and $S_9$; cuspidal specialization

When computing defining equations for minimally ramified three point covers, one normally works one cover at a time, and the equations obtained suggest that there do not exist simple general formulas. However, as an exception to this principle, one has the following two parameter family, with $N > m$ positive integers. Throughout, we put $r = N - m$.

$$\begin{aligned} \Lambda_{N,m} &= (N,\ 2\,1^{N-2},\ mr) \\ f_{N,m}(t,x) &= m^m x^N - t(Nx-r)^m \\ D_{N,m}(t) &= (-1)^{(N+2m)(N+1)/2} N^N m^{Nm} r^{(N-1)m} t^{N-1}(t-1). \end{aligned}$$

The cover $X_{N,m}$ is isomorphic to the cover $X_{N,r}$, an isomorphism being $x \mapsto mx/(Nx-r)$. So, without loss of generality, we restrict to the case $m \leq N/2$. We assume further that $N$ and $m$ are relatively prime; in this case the monodromy group $M_{N,m}$ and the Galois group $G_{N,m}$ are both the full symmetric group $S_N$.

For $m = 1$, the polynomial $f_{N,m}$ is a trinomial and these covers have been discussed in e.g. [Mat, Section II.3]. As we explain next, even the other covers $X_{N,m}$ can be given by trinomial equations, and in this guise they have appeared in many places. Define integers $v$ and $w$ by

$$\begin{aligned} v &\in \{0, 1, \ldots, r-1\} \\ Nv &\equiv 1 \bmod (r) \\ w &= (Nv-1)/r \end{aligned}$$

Define

$$y = \left(\frac{m}{Nx-r}\right)^{w-v} x^w \in \mathbf{Q}(x)$$

One has

$$y^m = t^{w-v} x \in \mathbf{Q}(x)$$

and $y$ is a root of

$$f^*_{N,m}(t,y) = my^N - Nt^v y^m + rt^w$$

In summary, we have two defining polynomials for the same cover, the canonical one $f_{N,m}$ and a trinomial version $f^*_{N,m}$. We use $f^*_{N,m}$ exclusively in the sequel. We mention $f_{N,m}$ because this is the cover given by the standard three-point cover algorithm sketched in Section 4: $x$ is a coordinate on $X$ while $y$ is a rational function in $x$, typically of degree $> 1$. The general trinomial $ax^N + bx^m + c$ fits in this situation via $t = (m/a)^m (b/N)^N (r/c)^r$. Thus we call the $X_{N,m}$ trinomial covers.

The bad reduction set $S_{N,m}$ is the set of primes dividing $Nmr$; in other words, $X_{N,m}$ has bad reduction within $S$ iff $N/m \in T_{\infty,\infty,\infty}(\mathbf{Z}^S)$. In particular, $S_{N,m} \subseteq \{2,3\}$ exactly for the four pairs $(N,m) = (2,1)$, $(3,1)$, $(4,1)$, and $(9,1)$. For comparison, we note that specializing $f_{3,1}$ just at $T^*_{\infty,2,\infty}$ already yields all 9 $A_3/S_3$ cubics. Similarly, specializing just at $T^*_{\infty,2,\infty}$ already yields all 23 $A_4/S_4$ quartics, except the unique totally real one; in fact, no specialization of $f_{4,1}$ is totally real.

We are interested here in the case $(N, m) = (9, 1)$, coming from the Catalan equation $2^3+1 = 3^2$. Here a Frobenius computation shows that nothing unexpected happens when one specializes $t$ to $\tau$ in the 35-element set $T^*_{8,2,9} = T^*_{\infty,2,\infty}$. One gets 10 fields with Galois group $A_9$ and 25 with Galois group $S_9$. The field absolute discriminants $2^a 3^b$ are mostly near the maximum possible $2^{31}3^{26}$, with $12 \leq a \leq 31$ and $12 \leq b \leq 26$ being the exact ranges. Formula 11.1 lets one understand $a$ and $b$ in terms of the 2-adic and 3-adic placement of $\tau$, respectively.

In general, suppose $X \to \mathbf{P}^1_{\mathbf{Q}}$ is a three point cover with defining equation $f(t, x)$. It makes sense to specialize $t$ to one of the cuspidal points $c = 0$, 1, and $\infty$ as well, to obtain a Galois field $K_c \subset \mathbf{C}$. Let $\mathrm{Norm}_c$ be the normalizer of the monodromy transformation $m_c$ in $G$. Then $\mathrm{Gal}(K_c/\mathbf{Q})$ is contained in $\mathrm{Norm}_c/\langle m_c\rangle$.

For the covers $X_{N,m}$ of this section, the most interesting cusp to specialize at is $c = 1$. The Galois group $G_{N,m,1}$ is contained in a symmetric group $S_{N-2}$. These fields $K_{N,m,1} \subset \mathbf{C}$ have been studied in [Bor]. In our example one has

$$f_{9,1}(1, x) = (x - 1)^2(x^7 + 2x^6 + 3x^5 + 4x^4 + 5x^3 + 6x^2 + 7x + 8).$$

The septic field defined by the degree seven factor has Galois group $S_7$ and discriminant $-2^{12}3^{10}$. We do not know whether the other nine $S_7$ fields on the complete list in [Jon] can be obtained from specializing three point covers.

Cuspidal specializations seem worthy of special attention. First, they depend only on the group-theoretic data defining the cover from which they come. Second, they serve as an aid in understanding non-cuspidal specializations according to the cuspidal control principle of Section 11. All the fields obtained from cuspidally specializing the covers in Sections 6-7 can be built from the low-degree fields tabulated in [JR1]. For example, $f_{27a}(1, x)$ and $f_{27c}(1, x)$ each factor as $f_6^2 f_{15}$, with $f_6$ having Galois group all of $S_6$; the factor $f_{15}$ has the same Galois group $S_6$, corresponding to a different permutation representation of $\mathrm{Gal}(\overline{\mathbf{Q}}/\mathbf{Q})$ with the same kernel. On the other hand, a cuspidal specialization from the cover $X_{28}$ of the next section goes slightly beyond the complete tables: the polynomial $f_{28}(1, x)$ factors as $f_6^2 f_{16}$; $f_6$ has Galois group $S_6$, serving as a resolvent for $f_{16}$, which has Galois group $2^4.S_6$.

## 9. $Sp_6(2)$; computing lower degree resolvents

Here we use the following new cover with $M = G = Sp_6(2) = W(E_7)'$:

$$\begin{aligned}
\Lambda_{28} &= (12^2\,3\,1,\ 2^6\,1^{16},\ 9^3\,1)\\
f_{28}(t, x) &= 3^6\,(x^2 + 6\,x + 6)^{12}x^3(3\,x + 4) - t\,2^{18}\,(x^3 + 3\,x^2 - 3)^9\\
D_{28}(t) &= 2^{540}3^{457}t^{24}(t - 1)^6.
\end{aligned}$$

Here $x$, as a multivalued function of $t$, can be expressed in terms of hypergeometric functions with finite monodromy, as this cover corresponds to entry 58 of the Beukers-Heckmann list [BH, page 353].

We specialize to $\tau$ in in the 35-element set $T^*_{12,2,9} = T^*_{\infty,2,\infty}$. A Frobenius computation proves that for $\tau \neq -1$ the algebras $L_{28,\tau}$ all have Galois group all of $Sp_6(2)$. At $\tau = -1$, there is an apparent group-drop to $S_8 \cong SO^+_6(2) \subset Sp_6(2)$. The rest of this section explains how we produce a degree eight polynomial with the same splitting field as $f_{28}(-1, x)$, thereby proving that indeed $G_{28,-1} \cong S_8$.

In general, suppose given an irreducible degree 28 polynomial $f$ over $\mathbf{Q}$ with Galois group $S_8$ or $Sp_6(2)$. For each bijection

$$\begin{array}{rcl} L = (\text{Two element subsets of } \{0, \ldots, 7\}) & \to & \text{Complex roots of f} \\ \{i,j\} & \to & \alpha_{L(i,j)} \end{array}$$

one has an octic polynomial

$$g_L(x) = \prod_{i=1}^{8}(x - \sum_{j\neq i} \alpha_{L(i,j)})$$

Assuming $f$ is sufficiently generic, one gets $28!/8! \approx 7.5 \times 10^{24}$ distinct octic polynomials. If the Galois group is $S_8$, then exactly one of these has coefficients in $\mathbf{Q}$, it being the desired polynomial.

Suppose now that $f$ has exactly four real roots so that the desired octic has two real roots. Then it suffices to consider bijections $L$ where the involution $-1$ of $\mathbf{Z}/8$ goes over to complex conjugation. Then one gets $(4!24!!)/(2!3!2^3) \approx 4.9 \times 10^{11}$ octic polynomials. So still this method for finding a resolvent octic is impractical.

In our case, however, we have an action of $Sp_6(2) \subset S_{28}$ on the roots of $f_{28}(-1, x)$ and the desired $g_L$ is among the thirty-six $g_L$ coming from the thirty-six subgroups of $Sp_6(2)$ isomorphic to $S_8$. Using monodromy techniques as described in Section 7, we can identify these 36 copies of $S_8$. Looking among only these, we find a labeling $L$ giving a desired $g$ with quite large coefficients. Applying *polred*, gives

$$g_8(x) = x^8 + 4x^7 + 8x^6 + 16x^5 + 22x^4 + 20x^3 + 10x^2 - 8x - 10, \tag{9.1}$$

with field discriminant $-2^{17}3^9$.

To prove rigorously that $f_{28}(-1, x)$ and $g_8(x)$ have the same splitting field, we proceed again as in Section 7. Write

$$\begin{array}{rcl} g_8(x) & = & \prod_{\alpha\in A'} (x - \alpha) \\ g_{28}(x) & = & \prod_{\{\alpha_1,\alpha_2\}\subset A'} (x - \alpha_1 - \alpha_2) \\ f_{28}(-1, x) & = & 3^7 \prod_{\beta\in B} (x - \beta) \end{array}$$

We need to numerically find the right bijection $\sigma$ from $A = \mathrm{Sub}_2(A')$ to $B$, and then algebraically confirm that one indeed gets an isomorphism from $\mathbf{Q}[x]/g_{28}(x)$ to $\mathbf{Q}[y]/f_{28}(-1, y)$. This unique correct bijection is indicated on the following chart, with the root-labelling convention of Section 7:

| | 1 | $2a$ | $2b$ | $3a$ | $3b$ | $4a$ | $4b$ | 5 |
|---|---|---|---|---|---|---|---|---|
| 1 | | $11b$ | $11a$ | $13a$ | $13b$ | $10a$ | $10b$ | 2 |
| $2a$ | | | 15 | $5b$ | $14b$ | $12b$ | $7b$ | $16a$ |
| $2b$ | | | | $14a$ | $5a$ | $7a$ | $12a$ | $16b$ |
| $3a$ | | | | | 9 | $4a$ | $8a$ | $3b$ |
| $3b$ | | | | | | $8b$ | $4b$ | $3a$ |
| $4a$ | | | | | | | 1 | $6a$ |
| $4b$ | | | | | | | | $6b$ |
| 5 | | | | | | | | |

Resolvents 7.1, 7.2, 7.3, 7.4 were computed and confirmed by this method.

## 10. $S_{32}$; PRIME-DROPPING SPECIALIZATION

Let $f(t,x) \in \mathbf{Q}[t,x]$ define a degree $N$ three-point cover, with bad reduction set $S'$. For $\tau \in \mathbf{Q} - \{0,1\}$ it may happen that the specialized algebra $L_\tau$ is ramified strictly within $S'$. Then we call $\tau$ a prime-dropping specialization point for $f$.

Consider again the trinomial covers $X_{N,m}$ of Section 8, with $r := N - m$. Recall that the specialized algebra is given by a trinomial equation: $L_\tau = \mathbf{Q}[x]/f^*_{N,m}(\tau, x)$. Recall also that the cover has bad reduction at all primes dividing $Nmr$. Nonetheless, one has the following fact, $i$ being any positive integer:

$$\text{(10.1)} \qquad \text{If } p^e || \left\{ \begin{array}{l} N \\ m \\ r \end{array} \right. \text{ and } \tau \in \left\{ \begin{array}{l} T(\mathbf{Q}_p)^{0,Nei} \\ T(\mathbf{Q}_p)^{\infty,mei} \\ T(\mathbf{Q}_p)^{\infty,rei} \end{array} \right. \text{ then } K_{N,m,\tau} \text{ is unramified at } p$$

This fact can be proved directly: for suitable rational numbers $af^*_{N,m}(\tau, bx) \in \mathbf{Z}[x]$, with polynomial discriminant prime to $p$, as in 10.2, 10.3, and 10.4 below. Alternatively, this fact is a special case of the main theorem of [LNV].

We found three new fields this way with $S = \{2,3\}$. All three fields come from the specialization point $\tau = 2 \cdot 5^5/3^2$, associated to the ABC triple

$$-2 \cdot 5^5 + 79^2 + 3^2 = 0.$$

The fields come from the polynomials

$$\begin{aligned} f^*_{8,3}(t,x) &= 3x^8 - 8t^2x^3 + 5t^3 \\ f^*_{9,4}(t,x) &= 4x^9 - 9t^4x^4 + 5t^7 \\ f^*_{32,5}(t,x) &= 5x^{32} - 32t^{11}x^5 + 27t^{13}. \end{aligned}$$

Nice defining equations for these number fields are as follows.

$$\text{(10.2)} \qquad 3^7 5^{-16} f^*_{8,3}(2 \cdot 5^5/3^2, -5^2x/3) = x^8 + 2^5x^3 + 2^3 3$$

$$\text{(10.3)} \qquad 2^{-2} 3^{18} 5^{-36} f^*_{9,4}(2 \cdot 5^5/3^2, 5^4x/3^2) = x^9 - 2^2 3^4 x^4 + 2^5 3^4$$

$$\text{(10.4)} \qquad 3^{32} 5^{65} f^*_{32,5}(2 \cdot 5^5/3^2, -5^2x/3) = x^{32} + 2^{16} 3^5 x^5 + 2^{13} 3^9.$$

The discriminants of these polynomials are $-2^{45}3^5 79^2$, $-2^{40}3^{48}79^2$, and $-2^{563}3^{277}79^2$, respectively. The field discriminants are $-2^{31}3^5$, $-2^{14}3^{24}$, and $-2^{191}3^{111}$, respectively. Here, in the degree 32 case, we used 11.1 to compute the exponents 191 and 111. Note that the maximal absolute discriminant allowed by local considerations in the last case is $2^{191}3^{112}$.

In all three cases, the fact that 79 does not divide the field discriminant $d$ can be seen directly from the polynomial, rather than by our usual appeal to 2.4. Namely, each polynomial factors in the form $f_{N-2}f_1^2$ over $\mathbf{F}_{79}$, with $f_{N-2}$ separable. This implies that $\mathrm{ord}_{79}(d) \in \{0,1\}$. However from the polynomial discriminant we know that $\mathrm{ord}_{79}(d) \in \{0,2\}$.

The method of constructing number fields ramified within $S$ by judiciously specializing covers with larger ramifications sets $S'$ gives many interesting fields in low degrees. However it does not seem too promising in general. It seems that for "the general" such cover $X$, all primes $p \in S'$ ramify in all specializations $L_\tau$. In other words zeroes, like the ones appearing on the $p = 2,3,5$ blocks of Table 11.1 below, seem to occur only exceptionally.

## 11. Variation of discriminants

Let $X$ be a three point cover, with defining polynomial $f(t,x)$. For $\tau \in \mathbf{Q}-\{0,1\}$ let $d(\tau)$ be the algebra discriminant of the specialized algebra $L_\tau$. An interesting question is how the $d(\tau)$ varies with $\tau$. This type of question has been emphasized by Serre, [Ser3, Section 6], [SD, Section 8.4.1]. Since the question is purely local, we replace $\mathbf{Q}$ by $\mathbf{Q}_p$ everywhere.

For trinomial covers, the main theorem of [LNV] goes a long way towards answering this question. Here is a somewhat weakened version of this theorem; our statement looks very different from the statement given in [LNV], because we make use of the coordinate $t$. *Fix relatively prime integers $N > m$ and put $r = N-m$. For $\tau \in T(\mathbf{Q}_p) = \mathbf{Q}_p - \{0,1\}$, let $d_{N,m}(\tau)$ be the algebra discriminant of $L_{N,m,\tau} = \mathbf{Q}_p[x]/f_{N,m}(\tau,x)$. Let $e = ord_p(Nmr) \geq 0$. Put $v = ord_p(\tau) \in \mathbf{Z}$. Then $ord_p(d_{N,m}(\tau))$ is given in many cases as follows.*

(11.1)

*Case $p^e||N$*

| | | |
|---|---|---|
| $m-(m,v-Ne)+r-(r,v-Ne)$ | *if* $v \in (-\infty,-Ne)$ | $(\tau \in T(\mathbf{Q}_p)^{\infty,-v})$ |
| $N(e+1)-(N,v)+v$ | *if* $v \in [-Ne,-1]'$ | $(\tau \in T(\mathbf{Q}_p)^{\infty,-v})$ |
| $N(e+1)-(N,v)$ | *if* $v \in [1,\infty)'$ | $(\tau \in T(\mathbf{Q}_p)^{0,v})$ |

*Case $p^e||m$*

| | | |
|---|---|---|
| $m(e+1)-(m,v)+r-(r,me-v)$ | *if* $v \in (-\infty,-1]'$ | $(\tau \in T(\mathbf{Q}_p)^{\infty,-v})$ |
| $m(e+1)-(m,v)+r-(r,me-v)-v$ | *if* $v \in [1,me]'$ | $(\tau \in T(\mathbf{Q}_p)^{0,v})$ |
| $N-(N,v-me)$ | *if* $v \in (me,\infty)$ | $(\tau \in T(\mathbf{Q}_p)^{0,v})$ |

*Case $p^e||r$*

| | | |
|---|---|---|
| $r(e+1)-(r,v)+m-(m,re-v)$ | *if* $v \in (-\infty,-1]'$ | $(\tau \in T(\mathbf{Q}_p)^{\infty,-v})$ |
| $r(e+1)-(r,v)+m-(m,re-v)-v$ | *if* $v \in [1,re]'$ | $(\tau \in T(\mathbf{Q}_p)^{0,v})$ |
| $N-(N,v-re)$ | *if* $v \in (re,\infty)$ | $(\tau \in T(\mathbf{Q}_p)^{0,v})$ |

*Here $(a,b)$ means $GCD(a,b)$. Also the symbol $'$ on an integer-interval indicates that integers divisible by $p$ are deleted.* Note that when $e=0$ the three formulas agree on their overlap, and together treat all cases with $v \neq 0$; moreover for $v \neq 0$ they agree with the general tame formula 2.4.

For arbitrary covers (not necessarily ramified only at three points) one has some interesting easily obtained information about the behavior of $L_\tau$ for $\tau \in \mathbf{P}^1(\mathbf{Q}_p)$ sufficiently near a cusp $c$. Suppose the partition of $N$ describing the local behavior of the cover $X$ near $c$ is $\prod m_k^k$, conventions here being that multiplicities here are labeled by bases and actual parts by exponents, so that e.g. $(N-2)^1 1^2$ corresponds to minimal ramification, like at $c=1$ in any $X_{N,m}$. Localize at the cusp $c$; at $c=0$ this means work over a ring of convergent power series $\mathbf{Q}_p\langle\langle t\rangle\rangle$, rather than over $\mathbf{Q}_p(t)$. The cover $X$ then splits into $\coprod X_k$, with $X_k$ having degree $km_k$ over the base curve (the $X_k$ themselves may well be disconnected, but we won't burden our notation by pursuing this point). So for $\tau$ sufficiently near $c$, one has a canonical decomposition $L_\tau = \prod L_{\tau,k}$. Also each $L_{\tau,k}$ contains a copy $L_{\tau,k,c}$ of the residual algebra $L_{c,k}$. Here $[L_{\tau,k} : L_{\tau,k,c}] = k$ and $[L_{\tau,k,c} : \mathbf{Q}_p] = m_k$. Finally there is a periodicity isomorphism, taking the form $L_{\tau,k} \cong L_{p^k\tau,k}$, for $c=0$. In brief, one has a cuspidal-control principle: *the behavior near a cusp is governed by the behavior at a cusp.*

Table 11.1 gives a more complete picture of the situation for the trinomial cover $X = X_{6,1}$. The left column gives values of $\mathrm{ord}_p(d(\tau))$ we have observed with $\tau \in T(\mathbf{Q}_p)^{\mathrm{gen}}$. Similarly, the printed number in the $(p, c, i)$ slot are values of $\mathrm{ord}_p(d(\tau))$ we have observed with $\tau \in T(\mathbf{Q}_p)^{c,i}$. For example, the *5, 3* in the left column means that in the generic 5-adic region 5 and 3 can arise as $\mathrm{ord}_5(d(\tau))$. Similarly, the *5, 3* in column 5 says that 5 and 3 have been observed to arise when $5^5$ exactly divides the denominator of $\tau$. The numbers in standard type are the only possibilities, according to Formula 11.1. For $p \geq 7$ the numbers in standard type are the only possibilities, ramification here being tame, and as described by 2.4. Since we have made a careful numerical investigation, we expect that the italicized numbers are the only possibilities as well. Also we expect that the periodic behavior has already begun where indicated by the parentheses. One conclusion from Table 11.1 is that the missing cases in Formula 11.1 are sometimes truly more complicated.

TABLE 11.1. Variation of field discriminants of $x^6 - 6tx + 5t$

| gen | $c$ | 1 | 2 | 3 | 4 | 5 | 6 | 7 | 8 | 9 | 10 | $\infty$ |
|---|---|---|---|---|---|---|---|---|---|---|---|---|
| $p = 2$ | 0 | 11 | *6, 8* | 9 | *6, 8* | (11 | *6, 0* | 11 | *4, 8* | 9 | *4, 8*) | $1^6$ |
| | 1 | (*7* | *6, 4*) | | | | | | | | | $1^2 4_4$ |
| | $\infty$ | 10 | *6, 4* | 6 | *6* | 6 | (0 | 4 | 4 | 4 | 4) | $1^5 1$ |
| $p = 3$ | 0 | 11 | 10 | (*7, 3* | 10 | 11 | *6, 2* | 11 | 10) | | | $1^6$ |
| | 1 | *7, 5* | *6* | (*3* | *4*) | | | | | | | $1^2 4_3$ |
| *6* | $\infty$ | 10 | 8 | *6* | 6 | 6 | (0 | 4 | 4 | 4 | 4) | $1^5 1$ |
| $p = 5$ | 0 | 8 | 7 | 6 | 5 | (0 | 5 | 4 | 3 | 4 | 5) | $1^6$ |
| | 1 | *6* | (*3* | *2*) | | | | | | | | $1^2 4_2$ |
| *5, 3* | $\infty$ | (9 | 9 | 9 | 9 | *5, 3*) | | | | | | $1^5 1$ |
| $p \geq 7$ | 0 | (5 | 4 | 3 | 4 | 5 | 0) | | | | | $1^6$ |
| | 1 | (1 | 0) | | | | | | | | | $1^2 4$ |
| 0 | $\infty$ | (4 | 4 | 4 | 4 | 0) | | | | | | $1^5 1$ |

The last column in Table 11.1 describes the cuspidal specializations $L_c$. In general, suppose $L$ is a separable $\mathbf{Q}_p$-algebra with degree $n$ and discriminant $p^d$, to be considered with multiplicity $k$. Then we indicate such an algebra by $n_d^k$, dropping subscripts $d = 0$ and superscripts $k = 1$. At the cusp $c = 1$, for example, one has the factorization

$$f_{6,1}(1, x) = (x - 1)^2(x^4 + 2x^3 + 3x^3 + 4x^4 + 5),$$

the quartic field having discriminant $2^4 3^3 5^2$. This accounts for our entry in the $c = 1$ rows. As $\tau$ goes to $c$, the minimum value of $d(\mathrm{ord}_p(L_\tau))$ is due to the $k = 1$ quartic while the oscillation above this minimal value is due to the $k = 2$ linear factor. In this way, each line of Table 11.1 gives an explicit example of the cuspidal control principle.

Since Katz motives behave so regularly in so many respects, we expect that ramification in corresponding three point covers $X$ may behave regularly too. If one had formulas generalizing 11.1, one would know the discriminants of the resulting algebras $L_\tau$ without computing defining equations. In favorable cases one might be able to deduce that $G_\tau$ is all of the generic group $G$ from the behavior of decomposition groups $D_p$, for $p$ running over the ramifying set $S$. As we said in our introduction, we think that this is a promising direction for further research.

## 12. Two general conjectures

We conclude with two general conjectures. We are confident that there some kernel of truth in both of them. However in the extremely strong form in which they are stated, they should be interpreted merely as the author's current best guess.

Let $\overline{\mathbf{Q}}_p$ be an algebraic closure of $\mathbf{Q}_p$. Let $\mathbf{Q}_p^{\mathrm{un}} \subset \overline{\mathbf{Q}}_p$ be the maximal unramified extension. Put $D_p = \mathrm{Gal}(\overline{\mathbf{Q}}_p/\mathbf{Q}_p)$ and $I_p = \mathrm{Gal}(\overline{\mathbf{Q}}_p/\mathbf{Q}_p^{\mathrm{un}})$. For $G$ a finite group let $\tilde{R}_p(G)$ be the set of homomorphisms from $I_p$ to $G$ which can be extended to $D_p$. The group $G$ acts on $\tilde{R}_p(G)$ by post-conjugation. Let $R_p(G)$ be the quotient; it is canonically defined, independent of the choice of $\overline{\mathbf{Q}}_p$. Put $r_p(G) = |R_p(G)|$.

In the tame case, i.e. when $p$ does not divide $|G|$, the set $R_p(G)$ can be explicitly described. Let $G^\natural$ be the set of conjugacy classes in $G$. Let $G^{\natural,p}$ be the subset consisting of classes $[\tau]$ with $[\tau] = [\tau^p]$. Then there is a canonical identification $R_p(G) = G^{\natural,p}$. In particular, for $p > n$, $R_p(S_n)$ can be identified with the set of partitions of $n$.

When $p$ does not divide $|G|$ one is in the wild case, and it is harder to describe $R_p(G)$ explicitly. However the sets $R_p(S_n)$ can be identified with the set of degree $n$ separable $p$-adic algebras, up to $\mathbf{Q}_p^{\mathrm{un}}$-isomorphism, and in [JR3] we deduce simple formulas for $r_p(S_n)$ from the Serre mass formula [Ser1].

Call a Galois extension of number fields $K_2/K_1$ Kummer if $\mathrm{Gal}(K_2/K_1)$ is abelian and $K_1$ contains all roots of unity dividing its exponent. Call a Galois number field $K \subset \mathbf{C}$ exceptional if there are subfields $K/K_2/K_1/\mathbf{Q}$ with $K_2/K_1$ Kummer with degree an odd prime. So the lowest-degree exceptional fields are $S_3$-fields containing $\mathbf{Q}(\sqrt{-3})$. Exceptional fields have to be treated separately as there are "too many of them"; we will simply exclude them here, writing $NF'(S,G)$ for the subset of $NF(S,G)$ consisting of non-exceptional fields. Note that for many pairs $(S,G)$, like all cases considered in this paper, $NF'(S,G)$ is all of $NF(S,G)$ for simple reasons.

The basic idea is then that there exists a constant $c_G$ such that for given $S$ one can expect

$$|NF'(S,G)| \approx c_G \prod_{p \in S} r_p(G).$$

One way to state this basic idea as a conjecture is as follows.

**Conjecture 12.1.** *Fix a finite group $G$ and let $S$ range over finite sets of primes where each prime $p$ in $S$ is required to be among the first $2|S|$ primes. Then the real numbers*

$$\frac{|NF'(S,G)|}{\prod_{p\in S} r_p(G)}$$

*have a unique limit point $c_G$.*

Let $J_p$ be the image of the inertia group $I_p$ in $\mathrm{Gal}(\overline{\mathbf{Q}}/\mathbf{Q})^{\mathrm{ab}}$. Then using the standard facts $\mathrm{Gal}(\overline{\mathbf{Q}}/\mathbf{Q})^{\mathrm{ab}} = \prod J_p$ and $J_p = \mathbf{Z}_p^\times$ it is easy to check that Conjecture 12.1 holds for for abelian $G$, with $c_G = |\mathrm{Aut}(G)|^{-1}$. In [JR3] we state refinements of this conjecture, and present numerical evidence for some solvable $G$ and some non-solvable small $G$. The aspect of Conjecture 12.1 in least doubt is the basic idea that $|NF(S,G)|$ grows quite rapidly with $S$. For example, from [Jon] we know that while $|NF(\{2,3\},S_5)| = 5$, $|NF(\{2,3,5\},S_5)| = 1353$.

We have not been able to formulate a general conjecture about the group constant $c_G$. We do expect however that $c_G$ should decrease with $|G|$; more precise statements might need to take into account the set of involutory conjugacy classes in $G$, these being the possible classes of complex conjugation.

In lieu of discussing $c_G$ further we fix $S = \{p,\ell\}$ and let $G$ vary. From Katz motives only, one can expect that e.g. $NF(\{p,\ell\}, GSp_n(\ell))$ is non-empty, for infinitely many $n$. Similarly, from classical modular forms of level a power of $p$ one can expect that e.g. $NF(\{p,\ell\}, PGL_n(\ell^e))$ is non-empty, for infinitely $e$. On the other hand, we think there is no analogous motivic or automorphic construction which naturally yields $A_n$ or $S_n$ fields for infinitely many $n$.

**Conjecture 12.2.** *For each finite set of primes $S$ there is a positive integer $n_S$ such that $NF(S,S_n) = \emptyset$ and $NF(S,A_n) = \emptyset$ for $n > n_S$.*

There all also natural constructions of number fields with controlled ramifying set which are neither motivic or automorphic in nature. These constructions tend to give preferred defining polynomials, explicit formulas for their polynomial discriminants, and $S_n$ or $A_n$ as Galois group. For example, one has classical polynomials coming from linear differential equations e.g. Jacobi polynomials and Bessel polynomials. More spectacularly, one has polynomials coming from non-linear differential equations, e.g. those in [NY]. In a different direction, one has moduli fields of three point covers, see e.g. [Mal3]. But all these constructions seem to yield only finitely many fields ramified within a given set $S$; this gives us some further confidence in Conjecture 12.2.

Department of Mathematics, Hill Center, Rutgers University, New Brunswick, NJ 08903

*E-mail address*: `davrobts@math.rutgers.edu`